\documentclass[11pt,a4paper]{amsart}
\usepackage{amsmath,amssymb,latexsym,amsthm,enumerate}
\usepackage{graphicx} 
\usepackage{hyperref} 
\usepackage{xcolor} 
\usepackage{psfrag} 

\theoremstyle{plain}
\newtheorem{theorem}{Theorem}[section]
\newtheorem{proposition}[theorem]{Proposition}
\newtheorem{lemma}[theorem]{Lemma}
\newtheorem{corollary}[theorem]{Corollary}
\theoremstyle{definition}
\newtheorem{definition}[theorem]{Definition}
\newtheorem{example}[theorem]{Example}

\theoremstyle{remark}
\newtheorem*{discussion}{Discussion}
\newtheorem*{impremark}{Important remark}
\newtheorem*{remark}{Remark}
\newtheorem*{remarks}{Remarks}

\newcounter{numpar}[section]

\newcommand*{\wt}{\widetilde}
\newcommand*{\ol}{\overline}

\newcommand*{\dd}{d}     

\newcommand*{\eps}{\varepsilon}

\newcommand*{\cE}{\mathcal E}
\newcommand*{\cF}{\mathcal F}

\newcommand*{\bbN}{\mathbb N}
\newcommand*{\bbQ}{\mathbb Q}
\newcommand*{\bbR}{\mathbb R}

\newcommand*{\fM}{\mathfrak M}

\newcommand*{\EE}{\mathsf E}
\newcommand*{\PP}{\mathsf P}
\newcommand*{\QQ}{\mathsf Q}

\newcommand*{\tPP}{{\wt\PP}}
\newcommand*{\tQQ}{{\wt\QQ}}

\DeclareMathOperator{\infline}{\ol\inf}

\newcommand*{\locll}{\mathrel{\stackrel{\mathrm{loc}}{\ll}}}
\newcommand*{\locsim}{\mathrel{\stackrel{\mathrm{loc}}{\sim}}}
\newcommand*{\LlocJ}{{L^1_{\mathrm{loc}}(J)}}
\newcommand*{\Llocr}{{L^1_{\mathrm{loc}}(r-)}}
\newcommand*{\Llocl}{{L^1_{\mathrm{loc}}(l+)}}
\newcommand*{\Llocinfty}{{L^1_{\mathrm{loc}}(\infty-)}}

\begin{document}
\title[On the Martingale Property]%
{On the Martingale Property of Certain Local Martingales}

\author{Aleksandar Mijatovi\'{c}}
\address{Department of Statistics, University of Warwick, UK}
\email{a.mijatovic@warwick.ac.uk}

\author{Mikhail Urusov}
\address{Institute of Mathematical Finance, Ulm University, Germany}
\email{mikhail.urusov@uni-ulm.de}

\thanks{We are grateful to Peter Bank, Nicholas Bingham, Mark Davis,
Yuri Kabanov, Ioannis Karatzas, Walter Schachermayer, and two anonymous referees
for valuable suggestions. This paper was written while the second author
was a postdoc in the Deutsche Bank Quantitative Products Laboratory, Berlin.}



\begin{abstract}
The stochastic exponential 
$Z_t=\exp\{M_t-M_0-(1/2)\langle M,M\rangle_t\}$
of a continuous local martingale 
$M$ 
is itself a continuous local martingale.
We give a necessary and sufficient condition for the process 
$Z$ 
to be a true martingale in the case 
where
$M_t=\int_0^t b(Y_u)\,dW_u$ 
and 
$Y$ 
is a one-dimensional diffusion driven by a Brownian motion~$W$.
Furthermore, we provide a necessary and sufficient condition
for 
$Z$ 
to be a uniformly integrable martingale 
in the same setting.
These conditions are deterministic
and expressed only in terms of the function 
$b$ 
and the drift and
diffusion coefficients 
of~$Y$.
As an application we provide a deterministic
criterion for the absence of bubbles in a one-dimensional setting.
\end{abstract}

\maketitle

\section{Introduction}
\label{I}
A random process that 
is a
local martingale 
but
does not satisfy the martingale property
is known as 
a \emph{strict local martingale}
(this terminology was introduced by
Elworthy, Li, and Yor~\cite{ElworthyLiYor:99}).
The question of whether a local martingale is
a strict local martingale
or a true martingale is of a particular interest for the stochastic exponential
\begin{equation}
\label{i0}
Z_t:=\cE(M)_t=\exp\left\{M_t-M_0-\frac12\langle M,M\rangle_t\right\}
\end{equation}
of a continuous local martingale $M$ because such a process $Z$ is often used
as a density process for a (locally) absolutely continuous measure change.
One can perform the measure change only if $Z$ is a martingale.
The problem of finding convenient
sufficient conditions on 
$M$ for $Z$  to be a martingale has attracted  significant interest
in the literature. 
The criteria of 
Novikov and Kazamaki are particularly well-known.
Novikov~\cite{Novikov:72} proved that the condition
\begin{equation}
\label{i2}
\EE\exp\left\{\frac12\langle M,M\rangle_t\right\}<\infty
\;\;\forall t\in[0,\infty)
\end{equation}
guarantees that $Z$ is a martingale.
Kazamaki~\cite{Kazamaki:77}
showed that $Z$ is a martingale provided
\begin{equation}
\label{i3}
\exp\left\{\frac12 M\right\}
\text{ is a submartingale.}
\end{equation}
Let us note that \eqref{i3} is equivalent to the condition
$\EE\exp\{(1/2)M_t\}<\infty$ for all $t\in[0,\infty)$
if $M$ is a true martingale (not just a continuous
local martingale as assumed above).
Novikov's criterion is of narrower scope
but often easier to apply.
For improvements on the criteria of  Novikov and Kazamaki in the setting
of Brownian motion
see e.g. Kramkov and Shiryaev~\cite{KramkovShiryaev:98},
Cherny and Shiryaev~\cite{ChernyShiryaev:01} and the references therein.
A similar question in the exponential semimartingale framework,
in particluar, for affine processes, has also attracted attention in the literature
(see e.g. Kallsen and Shiryaev~\cite{KallsenShiryaev:02},
Kallsen and Muhle-Karbe~\cite{KallsenMuhleKarbe:10},
Mayerhofer, Muhle-Karbe, and Smirnov~\cite{MayerhoferMuhleKarbeSmirnov:09}
and the references therein).
For treatments of related questions of (local) absolute continuity
of measures on filtered spaces see e.g. Jacod and Shiryaev~\cite[Ch. III and~IV]{JacodShiryaev:03}
and Cheridito, Filipovi\'c, and Yor~\cite{CheriditoFilipovicYor:05}.
While \eqref{i2} and \eqref{i3} are only sufficient conditions,
Engelbert and Senf~\cite{EngelbertSenf:91}
and, recently, Blei and Engelbert~\cite{BleiEngelbert:09}
provided necessary and sufficient conditions for $Z$ to be a martingale.
In~\cite{EngelbertSenf:91} the case of
a general continuous
local martingale $M$ 
is considered and the condition is given 
in terms of the time-change that turns 
$M$
into a (possibly stopped) Brownian motion.
In~\cite{BleiEngelbert:09} 
the case of a strong Markov
continuous local martingale~$M$
is studied
and the deterministic criterion
is expressed in terms of the speed measure of~$M$.
In the recent papers of Kotani~\cite{Kotani:06}
and Hulley and Platen~\cite{HulleyPlaten:08}
a related question is studied.
These authors obtain necessary and sufficient
conditions for a one-dimensional regular
strong Markov continuous local martingale to be a martingale.
For further literature review see e.g.
the bibliographical notes in the monographs
Karatzas and Shreve~\cite[Ch.~3]{KaratzasShreve:91},
Liptser and Shiryaev~\cite[Ch.~6]{LiptserShiryaev:01a},
Protter~\cite[Ch.~III]{Protter:05}, and
Revuz and Yor~\cite[Ch.~VIII]{RevuzYor:99}.

In the present paper we consider local martingales $M$
of the form  $M_t=\int_0^t b(Y_u)\,dW_u$, where $Y$ is a one-dimensional
diffusion driven by a Brownian motion~$W$.
Our main results are necessary and sufficient conditions
for $Z$ to be a true martingale (Theorem~\ref{MRA1})
and for $Z$ to be a uniformly integrable martingale (Theorem~\ref{MRA2}).
The conditions are deterministic and expressed only in terms of
the function $b$ and the drift and diffusion coefficients of~$Y$.

Compared with the aforementioned result of~\cite{EngelbertSenf:91},
our criterion is of narrower scope 
(Engelbert and Senf consider an arbitrary continuous local martingale~$M$),
but our results are easier to apply when 
$M_t=\int_0^t b(Y_u)\,dW_u$,
as the condition of 
Engelbert and Senf 
is given in terms of the Dambis--Dubins--Schwarz Brownian motion
of $M$ and the related time-change.
The setting of the present paper differs from that of~\cite{BleiEngelbert:09},
\cite{Kotani:06}, and~\cite{HulleyPlaten:08}
in that in our case
the process 
$Y$ 
possesses the
strong Markov property 
but both
$M=\int_0^\cdot b(Y_u)\,dW_u$
and 
$Z=\cE(M)$ 
can well be non-Markov.
As a simple example, consider $Y\equiv W$
and $b(x)=I(x>0)$. 
Then both 
$M$ 
and 
$Z$ 
have intervals
of constancy 
and therefore
the knowledge of the trajectory's past
helps to predict the future in the following way.
Let us fix a time $t$ and a position $M_t$ (resp.~$Z_t$).
If $M$ (resp.~$Z$) is constant immediately before~$t$,
$W$ has a negative excursion at time~$t$,
hence $M$ (resp.~$Z$) will be constant also immediately after~$t$.
Similarly, if $M$ or $Z$ is oscillating immediately before~$t$,
then it will be oscillating also immediately after~$t$.
This implies that $M$ and $Z$ are non-Markov.

We discuss the applications of Theorems \ref{MRA1} and~\ref{MRA2}
in specific situations by studying several examples.
There is evidence in the literature that in some settings
the loss of the martingale property of the stochastic exponential $Z$
is related to some auxiliary diffusion
exiting its state space 
(see Karatzas and Shreve~\cite[Ch.~5, Ex.~5.38]{KaratzasShreve:91},
Revuz and Yor~\cite[Ch.~IX, Ex.~(2.10)]{RevuzYor:99},
Sin~\cite{Sin:98}, Carr, Cherny, and Urusov~\cite{CarrChernyUrusov:07}).
Such a statement turns out to be true in the case where the diffusion
$Y$
does not exit its state space 
(see Corollary~\ref{MRA1.5})
but fails in general (see Example~\ref{ex:EAA1}).
The loss of the martingale property of $Z$ is in our setting related 
only to 
the auxiliary diffusion 
exiting its state space 
at a \emph{bad endpoint}
(this notion is defined in the next section).
As another application of 
Theorems \ref{MRA1} and~\ref{MRA2}, 
we obtain a deterministic necessary and sufficient condition
for the absence of bubbles in diffusion-based models.

The paper is structured as follows.
In Section~\ref{MR} we describe the setting and formulate the main results,
which are proved in Section~\ref{P}.
In Section~\ref{EAA} we illustrate the main results by a complete study of two examples.
Applications to financial bubbles are given in Section~\ref{MRB}.
In Section~\ref{ST} we recall the definition of \emph{separating time} (see~\cite{ChernyUrusov:04})
and formulate the results about separating times that are used in the proofs in Section~\ref{P}.
Finally, in the appendices we describe some facts about local times
and behaviour of solutions of SDEs with the coefficients satisfying the Engelbert--Schmidt conditions,
which are used in the paper.

\section{Main Results}
\label{MR}
We consider the state space $J=(l,r)$, $-\infty\le l<r\le\infty$
and a $J$-valued diffusion $Y=(Y_t)_{t\in[0,\infty)}$
on some filtered probability space
$(\Omega,\cF,(\cF_t)_{t\in[0,\infty)},\PP)$
governed by the SDE
\begin{equation}
\label{mr1}
dY_t=\mu(Y_t)\,dt+\sigma(Y_t)\,dW_t,\quad Y_0=x_0\in J,
\end{equation}
where $W$ is an $(\cF_t)$-Brownian motion and $\mu,\sigma\colon J\to\bbR$
are Borel functions  satisfying the Engelbert--Schmidt conditions
\begin{gather}
\label{mr2} 
\sigma(x)\ne0\;\;\forall x\in J,\\
\label{mr3}
\frac1{\sigma^2},\;\frac\mu{\sigma^2}\in\LlocJ.
\end{gather}
$\LlocJ$ denotes the class of locally integrable functions,
i.e. the functions $J\to\bbR$ that are integrable on compact subsets of~$J$.
Under conditions~\eqref{mr2} and~\eqref{mr3} SDE~\eqref{mr1}
has a unique in law  weak solution that possibly exits its 
state space $J$ 
(see~\cite{EngelbertSchmidt:85}, \cite{EngelbertSchmidt:91},
or \cite[Ch.~5, Th.~5.15]{KaratzasShreve:91}).
We denote the exit time of~$Y$
by $\zeta$. 
In the case 
$Y$ does exit its state space, i.e.
$\PP(\zeta<\infty)>0$, we need to specify
the behaviour of $Y$ after 
$\zeta$.
In what follows we assume that the solution
$Y$
on the set $\{\zeta<\infty\}$ 
stays after $\zeta$ at the boundary point of
$J$
at which it exits,
i.e. $l$ and $r$ become absorbing boundaries.
We will use the following terminology:
\begin{center}
\emph{$Y$ exits (the state space $J$) at $r$} means $\PP(\zeta<\infty,\lim_{t\uparrow\zeta}Y_t=r)>0$;

\emph{$Y$ exits (the state space $J$) at $l$} is understood in an analogous way.
\end{center}
\noindent
The Engelbert--Schmidt conditions are reasonable
weak assumptions. For instance, 
they are satisfied
if 
$\mu$ 
is locally bounded on 
$J$
and 
$\sigma$ 
is locally bounded away from zero on~$J$.
Finally, let us note that we assume neither that $(\cF_t)$ is generated by $W$
nor that $(\cF_t)$ is generated by~$Y$.

In this section we consider the stochastic exponential
\begin{equation}
\label{mra1}
Z_t=\exp\left\{\int_0^{t\wedge\zeta}b(Y_u)\,dW_u
-\frac12\int_0^{t\wedge\zeta}b^2(Y_u)\,du\right\},
\quad t\in[0,\infty),
\end{equation}
where we set $Z_t:=0$ for $t\ge\zeta$ on
$\{\zeta<\infty,\int_0^\zeta b^2(Y_u)\,du=\infty\}$.
In what follows we assume that
$b$ is a Borel function $J\to\bbR$ satisfying
\begin{equation}
\label{mra2}
\frac{b^2}{\sigma^2}\in\LlocJ.
\end{equation}
Below we show that condition~\eqref{mra2} is equivalent to
\begin{equation}
\label{mra3}
\int_0^t b^2(Y_u)\,du<\infty \quad\PP\text{-a.s. on }\{t<\zeta\},\quad t\in[0,\infty).
\end{equation}
Condition~\eqref{mra3} ensures that the stochastic integral
$\int_0^t b(Y_u)\,dW_u$ 
is well-defined on $\{t<\zeta\}$
and, as mentioned above,
is equivalent to imposing~\eqref{mra2} 
on the function~$b$.
Thus the process 
$Z=(Z_t)_{t\in[0,\infty)}$
defined in~\eqref{mra1}
is a nonnegative continuous local martingale
(continuity at time $\zeta$ on the set $\{\zeta<\infty,\int_0^\zeta b^2(Y_u)\,du=\infty\}$
follows from the Dambis--Dubins--Schwarz theorem; see~\cite[Ch.~V, Th.~1.6]{RevuzYor:99}).

\begin{discussion}
Since $Z$ is a nonnegative local martingale, it is a supermartingale
(by Fatou's lemma). Hence, for a fixed $T\in(0,\infty)$,
$Z$ is a martingale on the time interval $[0,T]$ if and only~if
\begin{equation}
\label{aa1}
\EE Z_T=1.
\end{equation}
As a nonnegative supermartingale, $Z$ has a ($\PP$-a.s.) limit
$Z_\infty:=\lim_{t\uparrow\infty}Z_t$ and is closed by $Z_\infty$,
i.e. the process $(Z_t)_{t\in[0,\infty]}$ (with time $\infty$ included)
is a supermartingale. Hence, $Z$ is a uniformly integrable martingale
if and only if \eqref{aa1} holds for $T=\infty$.
In Theorem~\ref{MRA1} below (see also the remark following the theorem)
we present a deterministic criterion
in terms of $\mu$, $\sigma$, and~$b$
for~\eqref{aa1} with $T\in(0,\infty)$.
In Theorem~\ref{MRA2} we give a deterministic
necessary and sufficient condition for~\eqref{aa1} with $T=\infty$.
\end{discussion}

Before we formulate the results let us show that \eqref{mra2} is equivalent to~\eqref{mra3}.
Since the diffusion 
$Y$
is a continuous semimartingale on the stochastic interval
$[0,\zeta)$,
by the occupation times formula 
we have
\begin{align}
\int_0^{t}b^2(Y_u)\,du
&=\int_0^{t}\frac{b^2}{\sigma^2}(Y_u)\,d\langle Y,Y\rangle_u\label{mra4}\\
&=\int_J \frac{b^2}{\sigma^2}(y)L_{t}^y(Y)\,dy,
\quad t\in[0,\zeta).\notag
\end{align}
The map
$(y,t)\mapsto L_t^y(Y)$,
defined on
$J\times[0,\zeta)$,
denotes a jointly continuous 
version of local time 
given in Proposition~\ref{LocTime_cont}.
For the definition and properties of the continuous semimartingale
local time see e.g.~\cite[Ch.~VI, \S~1]{RevuzYor:99}.
Now \eqref{mra3} follows from~\eqref{mra2} and the fact that
on $\{t<\zeta\}$ the function 
$y\mapsto L_{t}^y(Y)$
is  continuous 
and hence bounded
with a compact support in~$J$.
Conversely, suppose that \eqref{mra2} is not satisfied.
Then there exists a point $\alpha\in J$ such that
we either have
\begin{equation}
\label{mra5}
\int_\alpha^{\alpha+\eps}\frac{b^2}{\sigma^2}(y)\,dy=\infty\text{ for any }\eps>0
\end{equation}
or 
\begin{equation}
\label{mra6}
\int_{\alpha-\eps}^\alpha\frac{b^2}{\sigma^2}(y)\,dy=\infty\text{ for any }\eps>0.
\end{equation}
Below we assume \eqref{mra5} and $\alpha<x_0$ (recall that $x_0$ is the starting point for~$Y$).
Let us consider the stopping time
\begin{equation}
\label{mra7}
\tau_\alpha=\inf\{t\in[0,\infty)\colon Y_t=\alpha\}
\end{equation}
with the usual convention $\inf\emptyset:=\infty$.
By Proposition~\ref{APP0}, $\PP(\tau_\alpha<\zeta)>0$.
Hence, there exists $t\in(0,\infty)$ such that $\PP(\tau_\alpha<t<\zeta)>0$.
By Proposition~\ref{prop:a:lt2}, we have $L^y_t(Y)>0$ $\PP$-a.s. on the event $\{\tau_\alpha<t<\zeta\}$
for any $y\in[\alpha,x_0]$.
Since the mapping $y\mapsto L^y_t(Y)$ is continuous $\PP$-a.s. on $\{t<\zeta\}$
(see Proposition~\ref{LocTime_cont}), it is bounded away from zero on $[\alpha,x_0]$
$\PP$-a.s. on $\{\tau_\alpha<t<\zeta\}$.
Now it follows from \eqref{mra4} and~\eqref{mra5} that
$\int_0^t b^2(Y_u)\,du=\infty$ $\PP$-a.s. on $\{\tau_\alpha<t<\zeta\}$,
which implies that~\eqref{mra3} is not satisfied.
Other possibilities (\eqref{mra5} and $\alpha\ge x_0$,
\eqref{mra6} and $\alpha< x_0$, \eqref{mra6} and $\alpha\ge x_0$)
can be dealt with in a similar way.

Let us now consider an auxiliary $J$-valued diffusion $\wt Y$ governed by the SDE
\begin{equation}
\label{mra16}
d\wt Y_t=(\mu+b\sigma)(\wt Y_t)\,dt+\sigma(\wt Y_t)\,d\wt W_t,\quad \wt Y_0=x_0,
\end{equation}
on some probability space
$(\wt\Omega,\wt\cF,(\wt\cF_t)_{t\in[0,\infty)},\tPP)$.
SDE~\eqref{mra16} has a unique in law 
weak solution that possibly exits its state space, 
because the Engelbert--Schmidt conditions
\eqref{mr2} and~\eqref{mr3} are satisfied
for the coefficients $\mu+b\sigma$ and $\sigma$
(note that $b/\sigma\in\LlocJ$ due to~\eqref{mra2}).
Similarly to the case of the solution of SDE~\eqref{mr1} 
we denote the exit time of $\wt Y$
by $\wt\zeta$ and apply the following convention: 
on the set
$\{\wt\zeta<\infty\}$ 
the solution 
$\wt Y$
stays after $\wt\zeta$ at the boundary point at which it 
exits.

Let $\ol J:=[l,r]$.
Let us fix an arbitrary $c\in J$ and set
\begin{align}
\label{mra8}
\rho(x)&:=\exp\left\{-\int_c^x\frac{2\mu}{\sigma^2}(y)\,dy\right\},
\quad x\in J,\\
\label{mra10}
\wt\rho(x)&:=\rho(x)\exp\left\{-\int_c^x\frac{2b}{\sigma}(y)\,dy\right\},
\quad x\in J,\\
\label{mra9}
s(x)&:=\int_c^x\rho(y)\,dy,\quad x\in\ol J,\\
\label{mra11}
\wt s(x)&:=\int_c^x\wt\rho(y)\,dy,\quad x\in\ol J.
\end{align}
Note that $s$ (resp. $\wt s$) is the scale function
of diffusion~\eqref{mr1} (resp.~\eqref{mra16}).
Further, we set
\begin{equation*}
\wt v(x):=\int_c^x \frac{\wt s(x)-\wt s(y)}{\wt\rho(y)\sigma^2(y)}\,dy\quad\text{ and } \quad
v(x):=\int_c^x \frac{s(x)-s(y)}{\rho(y)\sigma^2(y)}\,dy,\quad x\in J.
\end{equation*}
Note that the functions
$\wt v$ 
and
$v$
are decreasing on $(l,c)$
and increasing on $(c,r)$. 
Therefore the quantities
\begin{equation*}
\wt v(r):=\lim_{x\uparrow r}\wt v(x),\quad
v(r):=\lim_{x\uparrow r}v(x),
\quad
\wt v(l):=\lim_{x\downarrow l}\wt v(x),\quad
v(l):=\lim_{x\downarrow l}v(x)
\end{equation*}
are well-defined.
By $\Llocr$ we denote the class of Borel functions $f\colon J\to\bbR$
such that $\int_x^r |f(y)|\,dy<\infty$ for some $x\in J$.
Similarly we introduce the notation $\Llocl$.

Let us recall that 
the process
$\wt Y$ 
(resp. 
$Y$)
exits its state space 
at the boundary point 
$r$
if and only if
\begin{equation}
\label{mra20}
\wt v(r)<\infty \quad\text{ (resp. } v(r)<\infty\text{)}.
\end{equation}
This is Feller's test for explosions (see \cite[Ch.~5, Th.~5.29]{KaratzasShreve:91}).
Sometimes it can be easier to check the following
condition
\begin{equation}
\label{mra21}
\wt s(r)<\infty
\text{ and }
\frac{\wt s(r)-\wt s}{\wt\rho\sigma^2}\in\Llocr
\quad
\text{(resp. } s(r)<\infty
\text{ and }
\frac{s(r)-s}{\rho\sigma^2}\in\Llocr\text{)},
\end{equation}
equivalent to~\eqref{mra20} 
(see \cite[Sec.~4.1]{ChernyEngelbert:05}).
Similarly, 
$\wt Y$ 
(resp. 
$Y$)
exits its state space 
at the boundary point 
$l$ 
if and only if
\begin{equation}
\label{mra22}
\wt v(l)<\infty \quad\text{ (resp. } v(l)<\infty\text{)}.
\end{equation}
Sometimes it can be easier to check the equivalent condition
\begin{equation}
\label{mra23}
\wt s(l)>-\infty
\text{ and }
\frac{\wt s-\wt s(l)}{\wt\rho\sigma^2}\in\Llocl
\quad
\text{(resp. } 
s(l)>-\infty
\text{ and }
\frac{s-s(l)}{\rho\sigma^2}\in\Llocl\text{)}.
\end{equation}

We say that the endpoint $r$ of $J$ is
\emph{good}
if
\begin{equation}
\label{mra12}
s(r)<\infty\text{ and }
\frac{(s(r)-s)b^2}{\rho\sigma^2}\in\Llocr.
\end{equation}
Sometimes it can be easier to check the equivalent condition
\begin{equation}
\label{mra13}
\wt s(r)<\infty\text{ and }
\frac{(\wt s(r)-\wt s)b^2}{\wt\rho\sigma^2}\in\Llocr.
\end{equation}
We say that the endpoint $l$ of $J$ is \emph{good} if
\begin{equation}
\label{mra14}
s(l)>-\infty\text{ and }
\frac{(s-s(l))b^2}{\rho\sigma^2}\in\Llocl.
\end{equation}
Sometimes it can be easier to check the equivalent condition
\begin{equation}
\label{mra15}
\wt s(l)>-\infty\text{ and }
\frac{(\wt s-\wt s(l))b^2}{\wt\rho\sigma^2}\in\Llocl.
\end{equation}
The equivalence between \eqref{mra12} and~\eqref{mra13}
as well as between \eqref{mra14} and~\eqref{mra15} 
will follow from remark~(iii) after Theorem~\ref{A0}.
If $l$ or~$r$ is not good, we call it \emph{bad}.

\begin{remark}
Even though 
conditions \eqref{mra12} and~\eqref{mra13} are equivalent,
the inequalities
$s(r)<\infty$ and $\wt s(r)<\infty$
are not equivalent.
The same holds for conditions \eqref{mra14} and~\eqref{mra15}.
\end{remark}

\begin{impremark}
The notions of good and bad endpoints of $J$
are central to the theorems below.
To apply these theorems we need to check
whether the endpoints are good in concrete situations.
In Section~\ref{ST} (see remark~(iv) following Theorem~\ref{A0})
we show that the endpoint $r$ is bad
whenever one of the processes $Y$ and $\wt Y$ exits at $r$
and the other does not.
This is helpful because one can sometimes immediately see that,
for example,
$Y$ does not exit at $r$ while $\wt Y$ does.
In such a case one can conclude that
$r$ is bad without
having to check either~\eqref{mra12} or~\eqref{mra13}.
The same holds for the endpoint~$l$.
\end{impremark}

\begin{theorem}
\label{MRA1}
Let the functions $\mu$, $\sigma$, and $b$
satisfy conditions \eqref{mr2}, \eqref{mr3}
and~\eqref{mra2}, and $Y$ be a
solution of SDE~\eqref{mr1}
that possibly exits the state space 
$J=(l,r)$.
Then the process $Z$ given by~\eqref{mra1}
is a martingale
if and only if at least one of the conditions (a)--(b) below is satisfied
AND at least one of the conditions (c)--(d) below is satisfied:

(a) $\wt Y$ does not exit $J$ at~$r$, i.e. \eqref{mra20}
\textup{(}equivalently,~\eqref{mra21}\textup{)} is not satisfied;

(b) $r$ is good, i.e. \eqref{mra12}
\textup{(}equivalently,~\eqref{mra13}\textup{)} is satisfied;

(c) $\wt Y$ does not exit $J$ at~$l$, i.e. \eqref{mra22}
\textup{(}equivalently,~\eqref{mra23}\textup{)} is not satisfied;

(d) $l$ is good, i.e. \eqref{mra14}
\textup{(}equivalently,~\eqref{mra15}\textup{)} is satisfied.
\end{theorem}

\begin{remark}
The same condition is necessary and sufficient
for $Z$ to be a martingale on the time interval $[0,T]$
for any fixed $T\in(0,\infty)$
(see the proof of Theorem~\ref{MRA1}
in Section~\ref{P}).
\end{remark}

The case of the diffusion 
$Y$,
which does not exit its state space, 
is of particular interest. In this case
Theorem~\ref{MRA1} takes a simpler form. 

\begin{corollary}
\label{MRA1.5}
Assume that 
$Y$
does not exit its state space and
let the assumptions of Theorem~\ref{MRA1} be satisfied.
Then $Z$ is a martingale if and only if
$\wt Y$ does not exit its state space. 
\end{corollary}

This corollary follows immediately from Theorem~\ref{MRA1}
and the important remark preceding it.

\begin{theorem}
\label{MRA2}
Under the assumptions of Theorem~\ref{MRA1}
the process $Z$ is a uniformly integrable martingale
if and only if at least one of the conditions (A)--(D) below is satisfied:

(A) $b=0$ a.e. on $J$ with respect to the Lebesgue measure;

(B) $r$ is good and $\wt s(l)=-\infty$;

(C) $l$ is good and $\wt s(r)=\infty$;

(D) $l$ and $r$ are good.
\end{theorem}

\begin{remark}
Condition~(A) in Theorem~\ref{MRA2} cannot be omitted.
Indeed, if $J=\bbR$, $b\equiv0$ and $Y:=W$ is a Brownian motion,
then 
$Z\equiv1$ 
is a uniformly integrable martingale
but none of the conditions (B), (C), (D) hold
because
$s(-\infty)=-\infty$ and $s(\infty)=\infty$
(and hence neither endpoint is good).
\end{remark}

\section{Examples}
\label{EAA}
\begin{example}
\label{ex:EAA1}
Let us fix $\alpha>-1$ and consider the state space
$J=(-\infty,\infty)$ and a diffusion $Y$ governed by the SDE
\begin{equation}
\label{eaa1}
dY_t=|Y_t|^\alpha\,dt+dW_t,\quad Y_0=x_0.
\end{equation}
For each $\alpha>-1$, the coefficients of~\eqref{eaa1}
satisfy the Engelbert--Schmidt conditions (see \eqref{mr2}--\eqref{mr3}).
Hence, SDE~\eqref{eaa1} has a unique in law 
weak solution that possibly exits its state space. 
We are interested in the stochastic exponential
\begin{equation*}
Z_t=\exp\left\{\int_0^{t\wedge\zeta}Y_u\,dW_u
-\frac12\int_0^{t\wedge\zeta}Y^2_u\,du\right\},
\quad t\in[0,\infty),
\end{equation*}
where we set $Z_t:=0$ for $t\ge\zeta$ on
$\{\zeta<\infty,\int_0^\zeta Y^2_u\,du=\infty\}$
(like in Section~\ref{MR}, $\zeta$ denotes the exit time of~$Y$).

Let us now apply the results of Section~\ref{MR}
to study the martingale property of~$Z$ (i.e. we have
$\mu(x)=|x|^\alpha$, $\sigma(x)\equiv1$, and $b(x)=x$).
We get the following classification:
\begin{enumerate}
\item If $-1<\alpha\le1$, then $Z$ is a martingale (but not uniformly integrable);
\item If $1<\alpha\le3$, then $Z$ is a strict local martingale;
\item If $\alpha>3$, then $Z$ is a uniformly integrable martingale.
\end{enumerate}

Let us now outline some key steps in getting these results
(computations are omitted). The auxiliary diffusion $\wt Y$
is in our case given by the SDE
\begin{equation*}
d\wt Y_t=(|\wt Y_t|^\alpha+\wt Y_t)\,dt+d\wt W_t,\quad \wt Y_0=x_0.
\end{equation*}
We have:
\begin{enumerate}
\item $\wt Y$ does not exit at~$-\infty$ (for each $\alpha>-1$);
\item If $-1<\alpha\le1$, then $\wt Y$ does not exit at~$\infty$;
\item If $\alpha>1$, then $\wt Y$ exits at~$\infty$;
\item $s(-\infty)=-\infty$, hence, the endpoint $-\infty$ is bad (for each $\alpha>-1$);
\item The endpoint $\infty$ is good if and only if $\alpha>3$;
\item If $\alpha>3$, then $\wt s(-\infty)=-\infty$.\footnote{In fact we have
$\wt s(-\infty)=-\infty$ if and only if $\alpha\ge1$.}
\end{enumerate}
This is sufficient for the application of the results in Section~\ref{MR},
which yield the classification above.
\end{example}

\begin{remarks}
(i) It is not surprising to see a qualitative change
in the properties of $Z$ as $\alpha$ passes level~$1$
because a qualitative change in the properties of $Y$ occurs
($Y$ is does not exit its state space for $\alpha\le1$ and exits at $\infty$ for $\alpha>1$).
However, it is much more surprising to see a qualitative change
in the properties of $Z$ as $\alpha$ passes level~$3$,
while $Y$ undergoes no qualitative change.

One can understand what happens for $\alpha>3$ as follows.
Let us recall that $Z$ is a martingale if and only if
it does not lose mass at finite times,
and is a uniformly integrable martingale if and only if
it does not lose mass at time infinity (see~\eqref{aa1}).
Informally, when $\alpha>3$, the process $Y$ exits its state space
``so quickly'' 
that $Z$ does not lose mass.

\smallskip
(ii) There is some evidence in the literature
(see Karatzas and Shreve~\cite[Ch.~5, Ex.~5.38]{KaratzasShreve:91},
Revuz and Yor~\cite[Ch.~IX, Ex.~(2.10)]{RevuzYor:99},
Sin~\cite{Sin:98},
Carr, Cherny, and Urusov~\cite{CarrChernyUrusov:07})
that the loss of the martingale property of $Z$ in some settings
is related to the exit of the auxiliary diffusion~$\wt Y$
from its state space.
The example studied above shows that this is no longer true in our setting.
In the case $\alpha>3$ the auxiliary diffusion $\wt Y$ exits the state space,
while $Z$ is a uniformly integrable martingale.
As we can see, the loss of the martingale property of $Z$ is related
only to the exit of $\wt Y$ from its state space at a bad endpoint.

\smallskip
(iii) Wong and Heyde~\cite{WongHeyde:04} consider processes $X$ on a filtration
generated by a Brownian motion $B$ that are progressively measurable functionals of~$B$
and study the martingale property of stochastic exponentials of $\int_0^. X_u\,dB_u$.
Let us note that the process $b(Y)$ in our setting is a progressively measurable functional of $W$
whenever $Y$ is a strong solution of~\eqref{mr1};
however, this may be not the case when $Y$ is a weak solution of~\eqref{mr1}
(see~\cite[Ch.~IX, Def.~(1.5)]{RevuzYor:99} for the employed terminology).
We remark that the process $Y$ in the present example is a (possibly explosive)
strong solution of~\eqref{eaa1} for $\alpha\ge1$
because the coefficients of~\eqref{eaa1} are locally Lipschitz
(see~\cite[Ch.~IX, Th.~(2.1) and Ex.~(2.10)]{RevuzYor:99}).
Hence, the stochastic exponential $Z$ in the present example for $\alpha\ge1$
fits exactly into the setting of Corollary~2 in~\cite{WongHeyde:04},
which deals with the processes that are (possibly explosive) solutions of SDEs.
Corollary~2 in~\cite{WongHeyde:04} appears to imply that $Z$ is a strict local martingale
if and only if the auxiliary diffusion $\wt Y$ exits its state space.
Since for $\alpha>3$ we know that $Z$ is a uniformly integrable martingale,
additional assumptions are required in Corollary 2 of~\cite{WongHeyde:04}
for the result to hold as stated.

The origin of this mistake in~\cite{WongHeyde:04} is their Proposition~1 (see \cite[p.~657]{WongHeyde:04}).
Namely, the formulation of Proposition~1 in~\cite{WongHeyde:04} is ambiguous
and its proof fails after the words ``which demonstrates the existence\ldots''
(see~\cite[p.~658, line~3]{WongHeyde:04}).
\end{remarks}

\begin{example}[Measure transformation in a generalised CEV process]
\label{ex:EAA2}
We now turn to a generalised constant elasticity
of variance (CEV) process that satisfies the SDE
\begin{equation}
\label{ea1}
dY_t=\mu_0 Y_t^\alpha dt+\sigma_0 Y_t^\beta dW_t,\quad Y_0=x_0\in J:=(0,\infty),
\quad\alpha,\beta\in\bbR,\>\mu_0\in\bbR\setminus\{0\},\>\sigma_0>0
\end{equation}
under a probability measure~$\PP$.
Note that for the chosen domain $J$
and the specified parameter ranges the coefficients in~\eqref{ea1}
are locally Lipschitz,
and therefore equation~\eqref{ea1} 
has a pathwise unique strong solution up to exit time~$\zeta$
(see~\cite[Ch.~IX, Ex.~(2.10)]{RevuzYor:99}).
We are interested in the stochastic exponential
\begin{equation}
\label{eq:eaa2}
Z_t=\exp\left\{
-\frac{\mu_0}{\sigma_0}\int_0^{t\wedge\zeta}Y_u^{\alpha-\beta}\,dW_u
-\frac12\frac{\mu_0^2}{\sigma_0^2}\int_0^{t\wedge\zeta}Y_u^{2\alpha-2\beta}\,du
\right\},\quad t\in[0,\infty),
\end{equation}
where we set $Z_t:=0$ for $t\ge\zeta$ on
$\{\zeta<\infty,\int_0^\zeta Y_u^{2\alpha-2\beta}\,du=\infty\}$.
If $Z$ is a martingale (resp. uniformly integrable martingale),
then we can perform a locally (resp. globally) absolutely continuous measure change
using $Z$ as the density process, and under the new measure $\QQ$
the process $Y$ will satisfy the driftless equation
$dY_t=\sigma_0Y_t^\beta\,dB_t$ (with a $\QQ$-Brownian motion~$B$),
i.e. $Y$ will be a classical CEV process under~$\QQ$.

Below we apply the results of Section~\ref{MR} to study the martingale property of~$Z$.
We have $\mu(x)=\mu_0 x^\alpha$, $\sigma(x)=\sigma_0 x^\beta$, and $b(x)=-\mu_0 x^{\alpha-\beta}/\sigma_0$.
Let us note that the case $\mu_0=0$ is trivial and therefore excluded in~\eqref{ea1}.
The auxiliary diffusion $\wt Y$ of~\eqref{mra16} satisfies the driftless SDE
$d\wt Y_t=\sigma_0\wt Y_t^\beta\,d\wt W_t$.
Then it follows that 
$\wt\rho(x)=1$
and
$\wt s(x)=x$
for all
$x\in J$
(note that conditions~\eqref{mra20}--\eqref{mra15}
are unaffected if we add a constant to~$\wt s$).
Since
$\wt s(\infty)=\infty$
the  process
$\wt Y$
does not exit at 
$\infty$
and,
by~\eqref{mra13},
the boundary point
$\infty$
is bad.
At~$0$
we have
$\wt s(0)>-\infty$
and hence, 
by~\eqref{mra23}, the process
$\wt Y$
exits at~$0$
if and only if 
$\beta<1$.
Similarly, by~\eqref{mra15},
we find that $0$
is a good boundary point if and only if
$2\beta-\alpha<1$.

Theorem~\ref{MRA1} 
yields that the process
$Z$
is a strict local martingale if and only if 
$0$
is a bad point
and the process
$\wt Y$
exits at~$0$.
This is equivalent to
the conditions
$2\beta-\alpha\ge1$
and
$\beta<1$.
Theorem~\ref{MRA2}
implies that 
$Z$
is a uniformly integrable martingale 
if and only if 
$0$
is a good boundary point,
in other words when the inequality
$2\beta-\alpha<1$
holds.
This analysis is graphically
depicted in Figure~\ref{EAF}.

\begin{figure}[htb!]
\centering
\psfrag{b1}{$\beta$} 
\psfrag{a1}{$\alpha$}
\includegraphics[width=0.7\textwidth]{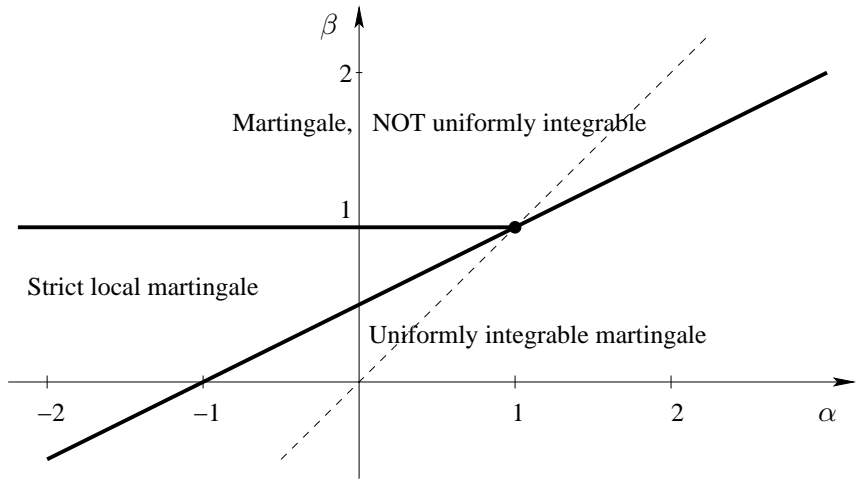}
\caption{\footnotesize The local martingale $Z$ defined in~\eqref{eq:eaa2}
is a uniformly integrable martingale
for all points in the open half-plane 
under the line
$2\beta-\alpha=1$,
a strict local martingale 
in the half-open wedge
defined by the inequalities
$\beta<1$
and
$2\beta-\alpha\ge1$,
and a martingale that is not uniformly integrable
in the closed wedge defined by the inequalities
$\beta\ge1$ and $2\beta-\alpha\ge1$.
Note that, 
excluding the trivial case
$\mu_0=0$,
the parameters
$\mu_0$
and 
$\sigma_0$
do not have any bearing on the 
martingale property of the process~$Z$.
The dot at 
$(1,1)$
corresponds to geometric Brownian motion and
the dashed line to the 
models with
$\alpha=\beta$.}
\label{EAF}
\end{figure}

In the case $\alpha=\beta$
the stochastic exponential~$Z$ defined in~\eqref{eq:eaa2}
is a geometric Brownian motion
stopped at the exit time $\zeta$ of the diffusion~$Y$.
As we can see in Figure~\ref{EAF},
$Z$ is a uniformly integrable martingale for $\alpha=\beta<1$
and a martingale that is not uniformly integrable for $\alpha=\beta\ge1$.
In particular the case
$\alpha=\beta=0$
implies that the law of Brownian motion with drift stopped at $0$
is absolutely continuous with respect to the law of Brownian motion 
stopped at $0$,
which is not the case if both processes are allowed to diffuse on the 
entire real line.
\end{example}

\section{Financial Bubbles}
\label{MRB}
Bubbles have recently attracted attention
in the mathematical finance literature;
see e.g. Cox and Hobson~\cite{CoxHobson:05},
Ekstr\"{o}m and Tysk~\cite{EkstromTysk:09},
Heston, Loewenstein, and Willard~\cite{HestonLoewensteinWillard:07},
Jarrow, Protter, and Shimbo~\cite{JarrowProtterShimbo:07},
Madan and Yor~\cite{MadanYor:06},
Pal and Protter~\cite{PalProtter:10}.
The reason is that option pricing in models
with bubbles is a very delicate task.
Therefore it is important to have tools
for ascertaining the existence or absence
of bubbles in financial models.
The problem of finding conditions that
guarantee the absence of bubbles in certain stochastic volatility models
was studied in Sin~\cite{Sin:98}, Jourdain~\cite{Jourdain:04},
Andersen and Piterbarg~\cite{AndersenPiterbarg:07};
see also the discussion in Lewis~\cite[Ch.~9]{Lewis:00}.
In this section we show how to apply Theorems \ref{MRA1} and~\ref{MRA2}
to get deterministic criteria for the absence of bubbles
in time-homogeneous local volatility models.

Let the discounted price of an asset be modelled by a nonnegative process~$S$.
Under a risk-neutral measure $\PP$ the discounted price is a local martingale.
In the terminology of~\cite{JarrowProtterShimbo:07}
there is a \emph{type~3 bubble} (resp. \emph{type~2 bubble})
in the model if $S$ is a strict local martingale
(resp. a martingale but not a uniformly integrable martingale) under~$\PP$.

\smallskip
The setting in this section is as follows.
Let $J=(0,\infty)$ and $Y$ be a $J$-valued diffusion
governed under a measure $\PP$ by the SDE
\begin{equation}
\label{mrba1}
dY_t=\mu_0 Y_t\,dt+\sigma(Y_t)\,dW_t,\quad Y_0=x_0>0,
\end{equation}
with $\mu_0\in\bbR$ and a Borel function $\sigma\colon J\to\bbR$
satisfying the Engelbert--Schmidt conditions
(as earlier $Y$ is stopped at the exit time~$\zeta$).
We interpret $Y$ as a non-discounted asset price
with evolution~\eqref{mrba1} under a risk-neutral measure.
Then $\mu_0$ is the risk-free interest rate (only
the case $\mu_0\ge0$ has a financial meaning,
but the results below hold for any $\mu_0\in\bbR$),
and the discounted price is $S_t=e^{-\mu_0 t}Y_t$.
Note that we also have $S_t=e^{-\mu_0 (t\wedge\zeta)}Y_t$
because $J=(0,\infty)$.

Since $Y$ is a price process, it is quite natural to assume that $Y$ does not exit at~$\infty$.
But this holds automatically for~\eqref{mrba1} as shown in the following lemma.
(This result is well-known for $\sigma$ having linear growth, but it holds also in general.
Surprisingly, we could not find this result in standard textbooks.
Let us also note that it is not at all easy to get this from Feller's test for explosions.)

\begin{lemma}
\label{MRBA1}
$Y$ does not exit 
the state space
$J=(0,\infty)$
at~$\infty$.
\end{lemma}

\begin{proof}
Until the exit time $\zeta$ we get by It\^o's formula
\begin{equation*}
e^{-\mu_0 t}Y_t=x_0+\int_0^t e^{-\mu_0 s}\sigma(Y_s)\,dW_s,\quad t\in[0,\zeta).
\end{equation*}
By the Dambis--Dubins--Schwarz theorem for continuous local martingales
on a stochastic interval (see \cite[Ch.~V, Ex.~(1.18)]{RevuzYor:99}),
the continuous local martingale $(e^{-\mu_0 t}Y_t)_{t\in[0,\zeta)}$
is a time-changed Brownian motion, hence
$\PP(\zeta<\infty,\lim_{t\uparrow\zeta}e^{-\mu_0 t}Y_t=\infty)=0$.
This concludes the proof.
\end{proof}

Note however that the process 
$Y$ 
may exit the interval
$J=(0,\infty)$
at~$0$.

\begin{proposition}
\label{MRBA2}
The discounted price process $S$ introduced above coincides with the process $x_0 Z$,
where $Z$ is given by~\eqref{mra1} with $b(x)=\sigma(x)/x$.
\end{proposition}

The desired criteria for the absence of type~2 and type~3 bubbles now follow immediately
from Proposition~\ref{MRBA2} and Theorems \ref{MRA1} and~\ref{MRA2}. We omit the formulations.

\begin{remark}
Note that in this setting formula~\eqref{mra10} simplifies to $\wt\rho(x)=\rho(x)/x^2$, $x\in J$,
since conditions \eqref{mra20}--\eqref{mra15} are not affected
if $\wt\rho$ is multiplied by a positive constant.
\end{remark}

\begin{proof}
Until the exit time $\zeta$ we set
\begin{equation*}
V_t:=\int_0^t \frac{dS_u}{S_u}=\int_0^t \frac{\sigma(Y_u)}{Y_u}\,dW_u,\quad t\in[0,\zeta).
\end{equation*}
The stochastic integral is well-defined because $Y$ is strictly positive on $[0,\zeta)$.
Then we have
\begin{equation*}
S_t=x_0 \cE(V)_t=x_0 Z_t,\quad t\in[0,\zeta).
\end{equation*}
Both $S$ and $Z$ have a limit as $t\uparrow\zeta$ and are stopped at~$\zeta$.
This completes the proof.
\end{proof}

In the important case of zero interest rate ($\mu_0=0$)
the problem under consideration amounts to the question of whether
the solution of the now driftless SDE 
in~\eqref{mrba1}
is a true martingale
or a strict local martingale.
This problem is of interest in itself,
and the answer looks particularly simple.
In what follows we will use, for a function $f\colon J\to\bbR$ and a class of functions~$\fM$,
the notation $f(x)\in\fM$ as a synonym for $f\in\fM$.

\begin{corollary}
\label{MRBA3}
Let $Y=(Y_t)_{t\in[0,\infty)}$
be a solution of the SDE $dY_t=\sigma(Y_t)\,dW_t$, $Y_0=x_0>0$,
where $\sigma\colon(0,\infty)\to\bbR$ is a Borel function satisfying
the Engelbert--Schmidt conditions
\textup{(}$Y$ is stopped at its hitting time of zero\textup{)}.

(i) If $x/\sigma^2(x)\in\Llocinfty$, then $Y$ is a strict local martingale.

(ii) If $x/\sigma^2(x)\notin\Llocinfty$, then $Y$ is a martingale but is not uniformly integrable.
\end{corollary}

\begin{remark}
Moreover, if $x/\sigma^2(x)\in\Llocinfty$,
then $Y$ is a strict local martingale
on each time interval $[0,T]$, $T\in(0,\infty)$.
\end{remark}

\begin{proof}
We can take $\rho\equiv1$, $s(x)=x$, $\wt\rho(x)=\frac1{x^2}$, $\wt s(x)=-\frac1x$.
Since $s(\infty)=\infty$ and $\wt s(0)=-\infty$,
neither $0$ nor $\infty$ are good (see \eqref{mra12} and~\eqref{mra15}),
the process $Y$ does not exit at $\infty$ and the diffusion $\wt Y$ does not exit at~$0$
(see~\eqref{mra21} and~\eqref{mra23}).
Theorem~\ref{MRA2} combined with Proposition~\ref{MRBA2}
yields that $Y$ cannot be a uniformly integrable martingale.
Theorem~\ref{MRA1} implies that $Y$ is a martingale
if and only if $\wt Y$ does not exit at~$\infty$.
The application of~\eqref{mra21} completes the proof.
Finally, the remark after Corollary~\ref{MRBA3}
follows from the remark after Theorem~\ref{MRA1}.
\end{proof}

The result established in Corollary~\ref{MRBA3}
was first proved in Delbaen and Shirakawa~\cite{DelbaenShirakawa:02b}
under the stronger assumption
that the functions $\sigma$ and $1/\sigma$ are locally bounded on $(0,\infty)$.
Later it was proved by a different method
in Carr, Cherny, and Urusov~\cite{CarrChernyUrusov:07}
without this assumption.
Furthermore, Corollary~\ref{MRBA3} can be obtained
as a consequence of Theorem~1 in Kotani~\cite{Kotani:06}
and of Theorem~3.20 in Hulley and Platen~\cite{HulleyPlaten:08}.
Let us note that the method of this paper
differs from all these approaches.

Let us further mention that it is only the possibility of bubbles in local volatility models
that prevents the main theorem in Mijatovi\'c~\cite{Mijatovic:10}
from being applied for the pricing of time-dependent barrier options.
In particular,
the above results could be used to characterise the class of local volatility models
within which time-dependent barrier options can be priced by solving a related system
of Volterra integral equations (see Section~2 in~\cite{Mijatovic:10}).

\begin{remarks}
(i) Let us discuss the question of whether all of the following possibilities
can be realised for the process $S_t=e^{-\mu_0 t}Y_t$, where $Y$ is given by~\eqref{mrba1}:

(1) $S$ is a strict local martingale;

(2) $S$ is a martingale but is not uniformly integrable;

(3) $S$ is a uniformly integrable martingale.

\noindent
Corollary~\ref{MRBA3} implies that the answer is affirmative
for (1) and~(2) even within the subclass $\mu_0=0$.
Possibility~(3) can also be realised:
one can take, for example, $\sigma(x)=\sqrt{x}$ and any $\mu_0>0$
(computations are omitted).

\smallskip
(ii) We now generalise Proposition~\ref{MRBA2}
to the case of a stochastic risk-free interest rate.
In the setting of Section~\ref{MR} let us assume that
$J=(l,r)$ with $0\le l<r\le\infty$,
so that the process $Y$ governed by~\eqref{mr1} is nonnegative.
Define the discounted price process until 
the exit time 
$\zeta$ 
by the formula
\begin{equation*}
S_t=\exp\left\{-\int_0^t \frac{\mu(Y_u)}{Y_u}\,du\right\}Y_t,
\quad t\in[0,\zeta).
\end{equation*}
By It\^{o}'s formula, $(S_t)_{t\in[0,\zeta)}$
is a nonnegative local martingale on the stochastic interval $[0,\zeta)$,
hence the limit $S_{\zeta}:=\lim_{t\uparrow\zeta}S_t$ exists
and is finite~$\PP$-a.s. (even if $r=\infty$ and $Y$ exits at~$\infty$).
Now we define a nonnegative local martingale $S=(S_t)_{t\in[0,\infty)}$
by stopping $(S_t)_{t\in[0,\zeta)}$ at $\zeta$ and interpret it as a discounted price
process.
In the case of a constant interest rate $\mu_0$ (i.e. $\mu(x)=\mu_0 x$)
and $J=(0,\infty)$ we get the setting above.
Proposition~\ref{MRBA2} holds without any change in this more general setting
and is proved in the same way.
\end{remarks}

\begin{example}[Bubbles in the CEV Model]
\label{ex:BCEV}
In the setting of this section let $\sigma(x)=\sigma_0 x^\alpha$,
$\sigma_0>0$, $\alpha\in\bbR$.
Then the process $Y$ given by~\eqref{mrba1}
is called the \emph{constant elasticity of variance (CEV) process}
(with drift~$\mu_0$).
This model was first proposed by Cox~\cite{Cox:75} for $\alpha\le1$
and Emanuel and MacBeth~\cite{EmanuelMacBeth:82} for $\alpha>1$
(see also Davydov and Linetsky~\cite{DavydovLinetsky:01}
and the references therein).
It is known that $S_t=e^{-\mu_0 t}Y_t$ is a true martingale
for $\alpha\le1$ and a strict local martingale for $\alpha>1$.
This was first proved in~\cite{EmanuelMacBeth:82}
by testing the equality $\EE S_t=x_0$
(the marginal densities of $Y$ are known explicitly
in this case). Below we prove this by our method
and further investigate when $S$ is a uniformly integrable
martingale.

Let us prove the following claims.
\begin{enumerate}[(i)]
\item $S$ is a martingale if and only if $\alpha\le1$;
\item $S$ is a uniformly integrable martingale
if and only if $\alpha<1$ and $\mu_0>0$.
\end{enumerate}

The case $\alpha=1$ is clear. We can therefore assume without loss of generality that $\alpha\ne1$.

By Proposition~\ref{MRBA2}
we need to study the martingale property
of the process $Z$ defined in~\eqref{mra1}
with $b(x):=\sigma_0 x^{\alpha-1}$, $x\in J:=(0,\infty)$.
The auxiliary diffusion in~\eqref{mra16}
takes the form
$$
d\wt Y_t=(\mu_0\wt Y_t+\sigma_0^2\wt Y_t^{2\alpha-1})\,dt
+\sigma_0\wt Y_t^\alpha\,d\wt W_t,\quad \wt Y_t=x_0.
$$
By Feller's test, $\wt Y$ exits at~$\infty$ if and only if $\alpha>1$.
It follows from the important remark preceding Theorem~\ref{MRA1} that
the endpoint $\infty$ is bad whenever $\wt Y$
exits at~$\infty$ (recall that $Y$ does not exit at~$\infty$).
By Theorem~\ref{MRA1}, $S$ is a strict local martingale for $\alpha>1$.

It remains to consider the case $\alpha<1$.
As was noted above, $\wt Y$ does not exit at $\infty$ in this case.
A direct computation shows that $\wt s(0)=-\infty$.
In particular, $\wt Y$ does not exit at~0.
By Theorem~\ref{MRA1} $S$ is a true martingale for $\alpha<1$.
Further, by Theorem~\ref{MRA2}, $S$ is a uniformly integrable
martingale if and only if $\infty$ is good.
Verifying~\eqref{mra12} (in the case $\alpha<1$)
we get that $\infty$ is good if and only if $\mu_0>0$.
This completes the proof of the claims (i) and~(ii).

\smallskip
It is interesting to note that the value of $\mu_0$
plays a role in the classification (i)--(ii) above.
Contrary to the case where $\mu_0=0$ (the CEV model without drift),
the possibility that $S$ is a uniformly integrable martingale
can be realised with positive interest rates.
In other words, taking 
$\mu_0>0$
may remove a type~2 bubble in this model,
but has no effect on the existence of
a type~3 bubble. 
\end{example}

\section{Separating Times}
\label{ST}
In order to present the proofs of our main results
we need the concept of a \emph{separating time}
for a pair of measures on a filtered space,
which was introduced in Cherny and Urusov~\cite{ChernyUrusov:04}.
In this section we recall the definition of separating times
and describe the explicit form of separating times
for distributions of solutions of SDEs.
All results stated in this section are taken from
Cherny and Urusov~\cite{ChernyUrusov:06}
(see also~\cite{ChernyUrusov:04}).

\subsection{Definition of Separating Times}
\label{STA}
Let $(\Omega,\cF)$ be a measurable space endowed with
a right-continuous filtration $(\cF_t)_{t\in[0,\infty)}$.
We recall that the $\sigma$-field $\cF_{\tau}$, for any $(\cF_t)$-stopping time~$\tau$, is defined by
\begin{equation}
\label{st0.5}
\cF_\tau=\{A\in\cF\colon A\cap\{\tau\le t\}\in\cF_t
\text{ for all }t\in[0,\infty)\}.
\end{equation}
In particular, $\cF_\infty=\cF$ by this definition.
Note that we do not assume here that 
$\cF=\bigvee_{t\in[0,\infty)}\cF_t$.
Definition~\eqref{st0.5}
is used for example in Jacod and Shiryaev~\cite{JacodShiryaev:03}
(we draw the reader's attention to the fact that there is
also an alternative definition of $\cF_{\tau}$ in the literature,
see e.g. Revuz and Yor~\cite{RevuzYor:99},
which is obtained by intersecting the $\sigma$-field in~\eqref{st0.5}
with the $\sigma$-field $\bigvee_{t\in[0,\infty)}\cF_t$).

Let $\PP$ and $\tPP$ be probability measures on~$(\Omega,\cF)$.
As usual, $\PP_\tau$ (resp.~$\tPP_\tau$)
denotes the restriction of $\PP$ (resp.~$\tPP$)
to~$(\Omega,\cF_{\tau})$.
In what follows, it will be convenient for us to consider
the extended positive half-line $[0,\infty]\cup\{\delta\}$, where
$\delta$ is an additional point.
We order $[0,\infty]\cup\{\delta\}$ in the following way:
we take the usual order on $[0,\infty]$ and let $\infty<\delta$.

\begin{definition}
\label{ST1}
An \emph{extended stopping time} is a map
$\tau\colon\Omega\to[0,\infty]\cup\{\delta\}$ such that
$\{\tau\le t\}\in\cF_t$ for any~$t\in[0,\infty]$.
\end{definition}

In order to introduce the notion of a separating time,
we need to formulate the following result.

\begin{proposition}
\label{ST2}
(i) There exists an extended stopping time~$S$ such
that, for any stopping time~$\tau$,
\begin{align}
\label{st1}
&\tPP_{\tau}\sim\PP_{\tau}\text{ on the set }\{\tau<S\}\,,\\
\label{st2}
&\tPP_{\tau}\perp\PP_{\tau}\text{ on the set }\{\tau\ge S\}\,.
\end{align}

(ii) If $S'$ is another extended stopping time with
these properties, then $S'=S$~$\PP,\tPP$-a.s.
\end{proposition}

\begin{definition}
\label{ST3}
A \emph{separating time} for $\PP$ and $\tPP$
(or, more precisely, for $(\Omega,\cF,(\cF_t),\PP,\tPP)$)
is an extended stopping time $S$ that satisfies~\eqref{st1}
and~\eqref{st2} for all stopping times~$\tau$.
\end{definition}

\begin{remark}
We stress that the separating time for $\PP$ and $\tPP$ is
determined $\PP,\tPP$-a.s.~uniquely.
\end{remark}

Informally, Proposition~\ref{ST2} states that two measures
$\PP$ and $\tPP$ are equivalent up to a random
time $S$ and become singular after that.
The equality $S=\delta$ means that $\PP$ and $\tPP$
never become singular, i.e. they are equivalent at time infinity.
On the contrary, $S=\infty$ means that $\PP$ and $\tPP$
are equivalent at finite times and become singular at time infinity.

Let us now provide a statistical interpretation of separating times,
which yields an intuitive understanding of the notion.
Suppose that we deal with the problem of sequentially distinguishing
between two statistical hypotheses $\PP$ and~$\tPP$,
where the information available to us at time $t$
is described by the $\sigma$-field~$\cF_t$.
(In particular, if the filtration $(\cF_t)_{t\in[0,\infty)}$
is the natural filtration of some process $X=(X_t)_{t\in[0,\infty)}$,
we observe sequentially a path of~$X$.)
Then, before the separating time $S$ occurs, we cannot 
know with certainty
what the true hypothesis is; and as soon as $S$ occurs we can
determine the true hypothesis with certainty.

In fact, the knowledge of the separating time for $\PP$ and $\tPP$
yields the knowledge of the mutual arrangement of $\PP$ and $\tPP$
from the viewpoint of their absolute continuity and singularity.
This is illustrated by the following result.
As usual, $\tPP\locll\PP$ (resp. $\tPP\locsim\PP$)
means that $\tPP_t\ll\PP_t$ (resp. $\tPP_t\sim\PP_t$)
for all $t\in[0,\infty)$.

\begin{lemma}
\label{ST4}
Let $S$ be a separating time for $\PP$ and~$\tPP$. Then

\noindent
\begin{equation*}
\begin{array}{llclclc}
{(i)}&\tPP\sim\PP&\Longleftrightarrow&S=\delta\;\;\PP,\tPP\text{-a.s.};&&&\\
{(ii)}&\tPP\ll\PP&\Longleftrightarrow&S=\delta\;\;\tPP\text{-a.s.};&&&\\
{(iii)}&\tPP\locsim\PP&\Longleftrightarrow&S\ge\infty\;\;\PP,\tPP\text{-a.s.};&&&\\
{(iv)}&\tPP\locll\PP&\Longleftrightarrow&S\ge\infty\;\;\tPP\text{-a.s.};&&&\\
{(v)}&\tPP\perp\PP&\Longleftrightarrow&S\le\infty\;\;\PP,\tPP\text{-a.s.}
&\Longleftrightarrow&S\le\infty\;\;\PP\text{-a.s.}&\phantom{cccccccccccccccccccccccccc}\\
{(vi)}&\tPP_0\!\perp\!\PP_0&\Longleftrightarrow&S=0\;\;\PP,\tPP\text{-a.s.}
&\Longleftrightarrow&S=0\;\;\PP\text{-a.s.}&\\
\end{array}
\end{equation*}
\end{lemma}

\begin{remark}
Other types of the mutual arrangement of $\PP$ and $\tPP$
are also easily expressed in terms of the separating time.
For example, we have
\begin{align*}
\tPP_t\ll\PP_t\;\;&\Longleftrightarrow\;\;S>t\;\;\tPP\text{-a.s.},\\
\tPP_t\perp\PP_t\;\;&\Longleftrightarrow\;\;S\le t\;\;\PP,\tPP\text{-a.s.}
\;\;\Longleftrightarrow\;\;S\le t\;\;\PP\text{-a.s.}
\end{align*}
for any~$t\in[0,\infty]$.
\end{remark}

The proof of Lemma~\ref{ST4} and the remark after it is straightforward.

\subsection{Separating Times for SDEs}
\label{STB}
Let us consider the state space $J=(l,r)$,
$-\infty\le l<r\le\infty$, and set $\ol J=[l,r]$.
In this subsection
we use the notation $(\Omega,\cF,(\cF_t)_{t\in[0,\infty)})$
for a certain canonical filtered space
and $\PP$ and $\tPP$ will be measures on this space,
namely, distributions of solutions of SDEs.
More precisely, let 
$\Omega:=\ol C([0,\infty),J)$ 
be the space
of continuous functions $\omega\colon[0,\infty)\to\ol J$
that start inside $J$ and can exit, i.e.
there exists $\zeta(\omega)\in(0,\infty]$
such that $\omega(t)\in J$ for $t<\zeta(\omega)$
and in the case $\zeta(\omega)<\infty$ we have
either $\omega(t)=r$ for $t\ge\zeta(\omega)$
(hence, also $\lim_{t\uparrow\zeta(\omega)}\omega(t)=r$)
or $\omega(t)=l$ for $t\ge\zeta(\omega)$
(hence, also $\lim_{t\uparrow\zeta(\omega)}\omega(t)=l$).
We denote the coordinate process on $\Omega$ by $X$
and consider the right-continuous canonical filtration
$\cF_t=\bigcap_{\eps>0}\sigma(X_s\colon s\in[0,t+\eps])$
and the $\sigma$-field $\cF=\bigvee_{t\in[0,\infty)}\cF_t$.
Note that the random variable $\zeta$ described above
is the exit time of~$X$.
Let the measures $\PP$ and $\tPP$ on $(\Omega,\cF)$
be the distributions of (unique in law) solutions of the SDEs
\begin{align}
dY_t&=\mu(Y_t)\,dt+\sigma(Y_t)\,dW_t,\quad Y_0=x_0,\label{stb1}\\
d\wt Y_t&=\wt\mu(\wt Y_t)\,dt+\wt\sigma(\wt Y_t)\,d\wt W_t,
\quad\wt Y_0=x_0,\label{stb2}
\end{align}
where the coefficients $\mu,\sigma$ as well as $\wt\mu,\wt\sigma$
are Borel functions $J\to\bbR$ satisfying
the Engelbert--Schmidt conditions \eqref{mr2}--\eqref{mr3}.

To summarize the setting, our input consists of the functions
$\mu$, $\sigma$, $\wt\mu$, and~$\wt\sigma$.
Given these functions we get the measures $\PP$ and $\tPP$
on our canonical space.
As an output we will present an explicit expression
for the separating time $S$ for $(\Omega,\cF,(\cF_t),\PP,\tPP)$.

Let $s$ denote the scale function of diffusion~\eqref{stb1}
and $\rho$ the derivative of~$s$ (see \eqref{mra8} and~\eqref{mra9}).
Similarly, let $\wt s$ be the scale function of~\eqref{stb2}
and $\wt\rho$ its derivative.
By $L^1_\mathrm{loc}(x)$ with $x\in J$ we denote
the set of Borel functions that are integrable
in a sufficiently small neighbourhood of~$x$.
Similarly we introduce the notation
$L^1_\mathrm{loc}(r-)$ and $L^1_\mathrm{loc}(l+)$.
Let $\nu_L$ denote the Lebesgue measure on~$J$.

We say that a point $x\in J$ is \emph{non-separating}
if $\sigma^2=\wt\sigma^2$ $\nu_L$-a.e. in a sufficiently small
neighbourhood of~$x$ and $(\mu-\wt\mu)^2/\sigma^4\in L^1_\mathrm{loc}(x)$.
We say that the right endpoint $r$ of $J$ is \emph{non-separating}
if all the points from $[x_0,r)$ are non-separating as well as
\begin{equation}
\label{stb3}
s(r)<\infty\text{ and }
(s(r)-s)\frac{(\mu-\wt\mu)^2}{\rho\sigma^4}\in L^1_\mathrm{loc}(r-).
\end{equation}
We say that the left endpoint $l$ of $J$ is \emph{non-separating}
if all the points from $(l,x_0]$ are non-separating as well as
\begin{equation}
\label{stb4}
s(l)>-\infty\text{ and }
(s-s(l))\frac{(\mu-\wt\mu)^2}{\rho\sigma^4}\in L^1_\mathrm{loc}(l+).
\end{equation}
A point in $\ol J$ that is not non-separating is called \emph{separating}.
Let $D$ denote the set of separating points in~$\ol J$.
Clearly, $D$ is closed in~$\ol J$.
Let us define
\begin{equation*}
D^\eps=\{x\in\ol J\colon\rho(x,D)<\eps\}
\end{equation*}
with the convention $\emptyset^\eps:=\emptyset$,
where $\rho(x,y)=|\arctan x-\arctan y|$, $x,y\in\ol J$
(the metric 
$\rho$ 
induces the standard topology on~$\ol J$;
we could not use here the standard Euclidean metric $|x-y|$,
as $\ol J$ can contain $\infty$ or~$-\infty$).

\begin{theorem}
\label{A0}
Let the functions $\mu,\sigma$ and $\wt\mu,\wt\sigma$
satisfy conditions \eqref{mr2}--\eqref{mr3}.
Let the measures $\PP$ and $\tPP$ on the canonical space $(\Omega,\cF)$
be the distributions of solutions of SDEs \eqref{stb1} and~\eqref{stb2}.
Then the separating time $S$ for $(\Omega,\cF,(\cF_t),\PP,\tPP)$
has the following form.

(i) If $\PP=\tPP$, then $S=\delta$~$\PP,\tPP$-a.s.

(ii) If $\PP\ne\tPP$, then
\begin{equation*}
S=\sup_n\infline\{t\in[0,\infty)\colon X_t\in D^{1/n}\}
\quad\PP,\tPP\text{-a.s.},
\end{equation*}
where ``$\,\infline$'' is the same as ``$\,\inf$''
except that $\infline\emptyset:=\delta$.
\end{theorem}

\begin{remarks}
(i) Theorem~\ref{A0} is proved in~\cite{ChernyUrusov:06}
in the particular case where $J=(-\infty,\infty)$
(see Theorem~5.7 in~\cite{ChernyUrusov:06}).
The case of general $J=(l,r)$ can be reduced
to the case $J=(-\infty,\infty)$
by considering a diffeomorphism $(l,r)\to(-\infty,\infty)$.
We omit the tedious but straightforward computations.

\smallskip
(ii) Let us explain the structure of $S$ in the case $\PP\ne\tPP$.
Denote by $\alpha$ the ``separating point that is
closest to $x_0$ from the left-hand side'', i.e.
\begin{equation*}
\alpha=\begin{cases}
\sup\{x:x\in[l,x_0]\cap D\}&\text{if}\;\;
[l,x_0]\cap D\ne\emptyset,\\
\Delta&\text{if}\;\;[l,x_0]\cap D=\emptyset,
\end{cases}
\end{equation*}
where $\Delta$ is an additional point ($\Delta\notin\ol J$).
Let us consider the ``hitting time of~$\alpha$'':
\begin{equation*}
U=\begin{cases}
\delta&\text{if}\;\;\alpha=\Delta,\\
\delta&\text{if}\;\;\alpha=l\text{ and }\liminf_{t\uparrow\zeta}X_t>l,\\
\zeta&\text{if}\;\;\alpha=l\text{ and }\liminf_{t\uparrow\zeta}X_t=l,\\
\ol\tau_\alpha&\text{if}\;\;\alpha>l,
\end{cases}
\end{equation*}
where $\ol\tau_\alpha=\infline\{t\in[0,\infty)\colon X_t=\alpha\}$.
Similarly, let us denote by $\beta$ the ``separating point that is
closest to $x_0$ from the right-hand side'' and by~$V$
the ``hitting time of~$\beta$''. Then $S=U\wedge V$~$\PP,\tPP$-a.s.
This follows from Proposition~\ref{APP1}.

\smallskip
(iii) We need the following statements proved in
remark~(ii) after Theorem~5.7 in~\cite{ChernyUrusov:06}.
If $[x_0,r)\subseteq\ol J\setminus D$, then
condition~\eqref{stb3} is equivalent to
\begin{equation}
\label{stb5}
\wt s(r)<\infty\text{ and }
(\wt s(r)-\wt s)\frac{(\mu-\wt\mu)^2}{\wt\rho\wt\sigma^4}\in L^1_\mathrm{loc}(r-).
\end{equation}
If $(l,x_0]\subseteq\ol J\setminus D$, then
condition~\eqref{stb4} is equivalent to
\begin{equation}
\label{stb6}
\wt s(l)>-\infty\text{ and }
(\wt s-\wt s(l))\frac{(\mu-\wt\mu)^2}{\wt\rho\wt\sigma^4}\in L^1_\mathrm{loc}(l+).
\end{equation}

In the context of Section~\ref{MR}
SDE~\eqref{stb2} has a particular form~\eqref{mra16},
i.e. we have $\wt\sigma=\sigma$ and $\wt\mu=\mu+b\sigma$.
Condition~\eqref{mra2} ensures that $J\subseteq\ol J\setminus D$.
Conditions \eqref{stb3}, \eqref{stb5}, \eqref{stb4}, and~\eqref{stb6}
reduce respectively to conditions
\eqref{mra12}, \eqref{mra13}, \eqref{mra14}, and~\eqref{mra15}.
We get the equivalence between \eqref{mra12} and~\eqref{mra13}
as well as between \eqref{mra14} and~\eqref{mra15}.

\smallskip
(iv) Let us show that the endpoint $r$ (resp.~$l$) is separating
whenever one of the diffusions \eqref{stb1} and~\eqref{stb2}
exits at $r$ (resp. at~$l$) and the other does not.
Suppose that $r$ is non-separating.
Since the set of non-separating points is open in~$\ol J$,
there exists $a<x_0$ such that all the points in $(a,r]$
are non-separating. Set
\begin{equation*}
E=\bigl\{\,\lim_{t\uparrow\zeta}X_t=r,\;X_t>a\;\forall t\in[0,\zeta)\bigr\}
\end{equation*}
and note that $\PP(E)>0$ by Proposition~\ref{APP2}.
By Theorem~\ref{A0}, we have $S=\delta$ on~$E$, hence $\tPP\sim\PP$ on~$E$.
If $X$ exits at $r$ under~$\PP$, then, by Proposition~\ref{APP2},
$\PP(E\cap\{\zeta<\infty\})>0$, hence $\tPP(E\cap\{\zeta<\infty\})>0$,
i.e. $X$ exits at $r$ also under~$\tPP$.
Similarly, if $X$ does not exit at $r$ under~$\PP$,
then it does not exit at $r$ also under~$\tPP$.

In the context of Section~\ref{MR}
(i.e. when $\wt\sigma=\sigma$ and $\wt\mu=\mu+b\sigma$)
we have that the endpoint $r$ is good (see~\eqref{mra12})
if and only if it is non-separating
(see~\eqref{stb3} and note that all points in $J$
are non-separating due to~\eqref{mra2}).
Thus, $r$ is bad whenever one of the diffusions \eqref{mr1} and~\eqref{mra16}
exits at $r$ and the other does not.
This proves the important remark preceding Theorem~\ref{MRA1}.
\end{remarks}

\section{Proofs of Theorems~\ref{MRA1} and~\ref{MRA2}}
\label{P}
In this section we prove the results formulated in Section~\ref{MR}.
We fix some number ${T\in(0,\infty)}$ and consider the following questions
for the process $Z$ defined in~\eqref{mra1}:

(i) Is the process $(Z_t)_{t\in[0,T]}$ a martingale?

(ii) Is the process $(Z_t)_{t\in[0,\infty)}$ a uniformly integrable martingale?

\noindent
The initial step is to translate these questions for the process $Z$ of~\eqref{mra1}
defined on an arbitrary filtered probability space into the related questions
on the canonical filtered space.
In contrast to Section~\ref{MR}, we use throughout this section
the notation $(\Omega,\cF,(\cF_t)_{t\in[0,\infty)})$ for the canonical
filtered space defined in Section~\ref{STB}.
By $X$ we denote the coordinate process on $\Omega$
and by $\zeta$ the exit time of~$X$.
Let the probability $\PP$ on $(\Omega,\cF)$
be the distribution of a (unique in law) solution of SDE~\eqref{mr1}.

Let us define a process $(W_t)_{t\in[0,\zeta)}$
on the stochastic interval $[0,\zeta)$
by the formula
\begin{equation*}
W_t=\int_0^t \frac1{\sigma(X_u)}\,(dX_u-\mu(X_u)\,du),
\quad t<\zeta.
\end{equation*}
Then it is a continuous $(\cF_t,\PP)$-local martingale
on the stochastic interval $[0,\zeta)$
with ${\langle W,W\rangle_t=t}$, $t<\zeta$.
Hence, $(W_t)_{t\in[0,\zeta)}$ can be extended
to a Brownian motion on the time interval $[0,\infty)$
on an enlargement of the probability space
(see \cite[Ch.~V, \S~1]{RevuzYor:99}).
Consequently, the process $(W_t)_{t\in[0,\zeta)}$ has
a finite limit $\PP$-a.s. on
$\{\zeta<\infty\}$
as $t\uparrow\zeta$,  
and we define
the process $W=(W_t)_{t\in[0,\infty)}$
by stopping $(W_t)_{t\in[0,\zeta)}$ at~$\zeta$
(and therefore do not enlarge the canonical probability space).
Thus, $W$ is a Brownian motion stopped at~$\zeta$.
Finally, we set
\begin{equation}
\label{p2}
Z_t=\exp\left\{\int_0^{t\wedge\zeta}b(X_u)\,dW_u
-\frac12\int_0^{t\wedge\zeta}b^2(X_u)\,du\right\},
\quad t\in[0,\infty)
\end{equation}
($Z_t:=0$ for $t\ge\zeta$ on $\{\zeta<\infty,\int_0^\zeta b^2(X_u)\,du=\infty\}$).
Thus the process $Z$ defined on $(\Omega,\cF,(\cF_t),\PP)$
is a nonnegative $(\cF_t,\PP)$-local martingale.

The key observation, which reduces the initial problem
to the canonical setting, is as follows:
the answer to the question~(i) (resp.~(ii))
for the process in~\eqref{mra1} is positive if and only if
the answer to the question~(i) (resp.~(ii))
for $Z$ of~\eqref{p2} is positive.
This may be difficult to believe at first sight since
the filtration in Section~\ref{MR}
need not be generated by~$Y$, while the canonical filtration $(\cF_t)$
in this section is the right-continuous filtration generated by~$X$.
However, this equivalence holds
because both in the case of~\eqref{mra1} and in the case of~\eqref{p2}
the property that $(Z_t)_{t\in[0,T]}$ is a martingale
(resp. $(Z_t)_{t\in[0,\infty)}$ is a uniformly integrable martingale)
is equivalent to the filtration-independent property $\EE Z_T=1$
(resp. $\EE Z_\infty=1$).
Thus, below we will consider the questions (i)--(ii)
for the process $Z$ of~\eqref{p2}.

Let $(a_n)_{n\in\bbN}$ and $(c_n)_{n\in\bbN}$ be strictly monotone sequences
such that $a_1<x_0<c_1$, $a_n\downarrow l$, and $c_n\uparrow r$. We set
\begin{equation*}
\tau_n=\inf\{t\in[0,\infty)\colon X_t\notin(a_n,c_n)\}
\end{equation*}
with $\inf\emptyset:=\infty$ and note that $\tau_n\uparrow\zeta$ for all $\omega\in\Omega$
and $\tau_n<\zeta$ on the set $\{\zeta<\infty\}$.
As in Section~\ref{ST} we denote by~$\cF_\tau$,
for any $(\cF_t)$-stopping time~$\tau$,
the $\sigma$-field defined in~\eqref{st0.5}.
Now we need to prove several lemmas.

\begin{lemma}
\label{PP1}
For any $t\in[0,\infty]$ we have
$\bigvee_{n\in\bbN}\cF_{t\wedge\tau_n}=\cF_t$.
\end{lemma}

Even though the statement in Lemma~\ref{PP1}
appears to be classical,
the fact that 
$\tau_n\uparrow\zeta$ 
may lead one to believe that it may fail to hold
on the event
$\{\zeta<t\}$.
The intuitive reason why the lemma holds is 
that the trajectory of the coordinate process $X$ up to $\zeta$
contains 
in itself (by the definition of $\Omega$)
the information about the whole trajectory of $X$.
The following proof formalises this intuitive idea.

\begin{proof}
We start by proving that
\begin{equation}
\label{pp1}
\bigvee_{n\in\bbN}\cF_{t\wedge\tau_n}=\cF\text{ on the set }\{\zeta\le t\}
\end{equation}
for any $t\in[0,\infty]$.
Indeed, this holds because $\cF=\sigma(X_s\colon s\in[0,\infty))$
and for any $s\in[0,\infty)$ we have
\begin{equation*}
X_s=\lim_{n\to\infty}X_{t\wedge\tau_n\wedge s}\text{ on }\{\zeta\le t\}.
\end{equation*}
The latter holds because $X$ stays after the exit time $\zeta$ at that endpoint of $J$
at which it exits (in particular, $X_s=X_{s\wedge\zeta}$).

The inclusion $\bigvee_{n\in\bbN}\cF_{t\wedge\tau_n}\subseteq\cF_t$ is clear.
For the reverse inclusion take an arbitrary set $A\in\cF_t$ and express it as
\begin{equation*}
A=\left[\bigcup_{n\in\bbN}(A\cap\{t<\tau_n\})\right]
\cup[A\cap\{\zeta\le t\}].
\end{equation*}
The set $A\cap\{t<\tau_n\}$ belongs both to $\cF_t$ and to $\cF_{\tau_n}$,
hence to $\cF_t\cap\cF_{\tau_n}=\cF_{t\wedge\tau_n}$.
Finally, by~\eqref{pp1},
$A\cap\{\zeta\le t\}\in\bigvee_{n\in\bbN}\cF_{t\wedge\tau_n}$.
\end{proof}

\begin{lemma}
\label{PP2}
For any $n\in\bbN$ we have
$\bigvee_{t\in[0,\infty)}\cF_{t\wedge\tau_n}=\cF_{\tau_n}$.
\end{lemma}

The proof is similar to that of Lemma~\ref{PP1}.

In what follows the probability measure $\tPP$ on $(\Omega,\cF)$
is the distribution of a (unique in law) solution of~\eqref{mra16}.
As in Section~\ref{ST}, for any $(\cF_t)$-stopping time~$\tau$,
$\PP_\tau$ (resp.~$\tPP_\tau$) denotes the restriction of~$\PP$
(resp.~$\tPP$) to the $\sigma$-field~$\cF_\tau$.
By $S$ we denote the separating time for $(\Omega,\cF,(\cF_t),\PP,\tPP)$
(see Definition~\ref{ST3}).

\begin{lemma}
\label{PP3}
For any $n\in\bbN$ we have $\tPP_{\tau_n}\sim\PP_{\tau_n}$ and
\begin{equation}
\label{p7}
\frac{d\tPP_{t\wedge\tau_n}}{d\PP_{t\wedge\tau_n}}
=Z_{t\wedge\tau_n}\quad\PP\text{-a.s.},\quad t\in[0,\infty].
\end{equation}
\end{lemma}

\begin{proof}
Let us fix $n\in\bbN$.
We want to apply Theorem~\ref{A0} to our measures $\PP$ and~$\tPP$.
In our case all points inside $J$ are non-separating due to~\eqref{mra2}.
Hence,
\begin{equation*}
S\ge\zeta>\tau_n\quad\PP,\tPP\text{-a.s.},
\end{equation*}
where the first inequality follows from remark~(ii) after Theorem~\ref{A0}
and the second one from Proposition~\ref{APP1}.
By the definition of a separating time, we get $\tPP_{\tau_n}\sim\PP_{\tau_n}$.
It remains to prove~\eqref{p7}.

We consider a c\`adl\`ag version of the density process
$D^{(n)}_t:=\frac{d\tPP_{t\wedge\tau_n}}{d\PP_{t\wedge\tau_n}}$, $t\in[0,\infty)$,
which is an $(\cF_{t\wedge\tau_n},\PP)$-martingale
closed by the $\PP$-a.s. strictly positive random variable~$\frac{d\tPP_{\tau_n}}{d\PP_{\tau_n}}$.
Then the processes $D^{(n)}$ and $D^{(n)}_-$ are $\PP$-a.s. strictly positive
(see \cite[Ch.~III, Lem.~3.6]{JacodShiryaev:03}).
Since $\frac{d\tPP_{\tau_n}}{d\PP_{\tau_n}}$ is $\cF_{\tau_n}$-measurable, we have
\begin{equation*}
D^{(n)}_t=\EE_\PP\bigg(\,\frac{d\tPP_{\tau_n}}{d\PP_{\tau_n}}\big|\cF_{t\wedge\tau_n}\bigg)
=\EE_\PP\bigg(\,\frac{d\tPP_{\tau_n}}{d\PP_{\tau_n}}\big|\cF_t\bigg)\quad\PP\text{-a.s.},
\quad t\in[0,\infty).
\end{equation*}
Hence, $D^{(n)}$ is also an $(\cF_t,\PP)$-martingale.
We obtain that the stochastic logarithm
\begin{equation}
\label{pp2}
L^{(n)}_t:=\int_0^t \frac{dD^{(n)}_u}{D^{(n)}_{u-}},\quad t\in[0,\infty),
\end{equation}
is a well-defined $(\cF_t,\PP)$-local martingale stopped at~$\tau_n$
(as $D^{(n)}$ is stopped at~$\tau_n$).
Let us prove that
\begin{equation}
\label{eq:FromL}
L^{(n)}_t=\int_0^{t\wedge\tau_n}b(X_u)\,dW_u\quad\PP\text{-a.s.},
\quad t\in[0,\infty).
\end{equation}

To prove~\eqref{eq:FromL} we first need to argue that there exists
a predictable process $H^{(n)}$
integrable with respect to the $(\cF_t,\PP)$-local martingale
\begin{equation*}
M^{(n)}_t:=X_{t\wedge\tau_n}-\int_0^{t\wedge\tau_n}\mu(X_u)\,du,\quad t\in[0,\infty),
\end{equation*}
such that
\begin{equation}
\label{p10}
L^{(n)}_t=\int_0^t H^{(n)}_u\,dM^{(n)}_u\quad\PP\text{-a.s.},
\quad t\in[0,\infty).
\end{equation}
Note that existence of such a process
$H^{(n)}$ will imply that both $L^{(n)}$ and $D^{(n)}$ are \mbox{$\PP$-a.s.} continuous processes
(as $D^{(n)}$ is the stochastic exponential of~$L^{(n)}$).
The idea is to use 
the Fundamental Representation Theorem
(see \cite[Ch.~III, Th.~4.29]{JacodShiryaev:03})
in order to prove the existence of~$H^{(n)}$.
However, we cannot apply the Fundamental Representation Theorem directly
because the process $X$ may exit the state space $J$ under~$\PP$
and thus may not be a semimartingale.

In order to avoid this problem
we consider a probability measure~$\PP'$
on the canonical filtered space,
which is the distribution of a (unique in law) solution of the SDE
\begin{equation}
\label{eq:AuxililarySDE}
dY_t=\mu'(Y_t)\,dt+\sigma'(Y_t)\,dB_t,\quad Y_0=x_0,
\end{equation}
where $B$ denotes some Brownian motion.
The Borel functions 
$\mu',\sigma'\colon J\to\bbR$
are given by
$\mu'=\mu$ and $\sigma'=\sigma$ on $[a_{n+1},c_{n+1}]$.
On the complement $J\setminus[a_{n+1},c_{n+1}]$
they are chosen so that
the Engelbert--Schmidt conditions \eqref{mr2} and~\eqref{mr3}
are satisfied and so that the coordinate process $X$
does not exit the state space $J$
under~$\PP'$.
Then by~\cite[Th.~2.11]{ChernyEngelbert:05}
we obtain the equality $\PP_{\tau_n}=\PP'_{\tau_n}$.
Since the process $L^{(n)}$ is stopped at~$\tau_n$,
it is $\cF_{\tau_n}$-measurable.
Let $\eta_m\uparrow\infty$~$\PP$-a.s. be a localizing sequence for~$L^{(n)}$,
i.e. $(L^{(n)})^{\eta_m}$ are $(\cF_t,\PP)$-martingales.
If we define
$\eta'_m:=\eta_m\wedge\tau_n$ and use the fact
$\PP_{\tau_n}=\PP'_{\tau_n}$,
we find that $\eta'_m\uparrow\tau_n$~$\PP'$-a.s. as $m\uparrow\infty$
and that $(L^{(n)})^{\eta'_m}$ is an $(\cF_t,\PP')$-martingale.
Therefore, $L^{(n)}$ is an $(\cF_t,\PP')$-local martingale
on the stochastic interval $[0,\tau_n)$
(we cannot guarantee more knowing only $\PP_{\tau_n}=\PP'_{\tau_n}$;
for instance, we cannot guarantee that $\eta_m\uparrow\infty$~$\PP'$-a.s.
because $\PP'$ need not be locally absolutely continuous with respect to~$\PP$).
Now the Fundamental Representation Theorem applied to all
$(L^{(n)})^{\eta'_m}$, $m\in\bbN$, implies the existence of
a predictable process $H'$ on the filtered probability space
$(\Omega,\cF,(\cF_t)_{t\in[0,\infty)},\PP')$
such that
\begin{equation}
\label{p12}
L^{(n)}_t=\int_0^t H_u'\,(dX_u-\mu'(X_u)\,du)\quad\PP'\text{-a.s.},
\quad t\in[0,\tau_n),
\end{equation}
since uniqueness in law holds for SDE~\eqref{eq:AuxililarySDE}.
We now obtain~\eqref{p10} by setting $H^{(n)}_u:=H'_u I(u\le\tau_n)$,
stopping the integrator in~\eqref{p12} at $\tau_n$,
and using the equality $\PP_{\tau_n}=\PP'_{\tau_n}$.

We get from~\eqref{pp2} that
\begin{equation}
\label{pp3}
\frac{d\tPP_{t\wedge\tau_n}}{d\PP_{t\wedge\tau_n}}
=\cE(L^{(n)})_t\quad\PP\text{-a.s.},\quad t\in[0,\infty).
\end{equation}
The process $M^{(n)}$ is a continuous $(\cF_t,\PP)$-local martingale stopped at~$\tau_n$
and is therefore also an $(\cF_{t\wedge\tau_n},\PP)$-local martingale.
By Girsanov's theorem for local martingales (see \cite[Ch.~III, Th.~3.11]{JacodShiryaev:03})
applied on the filtration $(\cF_{t\wedge\tau_n})_{t\in[0,\infty)}$
and by formulas \eqref{pp3} and~\eqref{p10}, the process
\begin{equation*}
M^{(n)}_t-\langle M^{(n)},L^{(n)}\rangle_t
=X_{t\wedge\tau_n}-\int_0^{t\wedge\tau_n}
(\mu(X_u)+H^{(n)}_u \sigma^2(X_u))\,du,\quad t\in[0,\infty),
\end{equation*}
is an $(\cF_{t\wedge\tau_n},\tPP)$-local martingale.
Since $\tPP$ is the distribution of a solution of SDE~\eqref{mra16},
the process
\begin{equation*}
X_{t\wedge\tau_n}-\int_0^{t\wedge\tau_n}(\mu+b\sigma)(X_u)\,du,\quad t\in[0,\infty),
\end{equation*}
is an $(\cF_t,\tPP)$-local martingale and therefore
also an $(\cF_{t\wedge\tau_n},\tPP)$-local martingale.
The process obtained as a difference of these two
$(\cF_{t\wedge\tau_n},\tPP)$-local martingales
\begin{equation*}
\int_0^{t\wedge\tau_n}
(b(X_u)\sigma(X_u)-H^{(n)}_u \sigma^2(X_u))\,du,
\quad t\in[0,\infty),
\end{equation*}
is a continuous $(\cF_{t\wedge\tau_n},\tPP)$-local martingale
of finite variation starting from zero.
Hence, it is identically zero.
Consequently, $\tPP$-a.s. we have
\begin{equation}
\label{p14}
H^{(n)}_u=\frac b\sigma(X_u)
\text{ for }\nu_L\text{-a.e. }u\in[0,\tau_n],
\end{equation}
where $\nu_L$ denotes the Lebesgue measure.
Since $\tPP_{\tau_n}\sim\PP_{\tau_n}$,
equality~\eqref{p14} holds also~$\PP$-a.s.,
and \eqref{eq:FromL} follows from~\eqref{p10}.
Now \eqref{p7} for $t\in[0,\infty)$ follows from \eqref{pp3} and~\eqref{eq:FromL},
and it remains only to prove \eqref{p7} for $t=\infty$.

By Levy's theorem and Lemma~\ref{PP2},
\begin{equation*}
\frac{d\tPP_{t\wedge\tau_n}}{d\PP_{t\wedge\tau_n}}\xrightarrow[t\uparrow\infty]{}
\EE_\PP\bigg(\,\frac{d\tPP_{\tau_n}}{d\PP_{\tau_n}}
\Big|\bigvee_{t\in[0,\infty)}\cF_{t\wedge\tau_n}\bigg)
=\frac{d\tPP_{\tau_n}}{d\PP_{\tau_n}}\quad\PP\text{-a.s.}
\end{equation*}
Thus, \eqref{p7} for $t\in[0,\infty)$ implies \eqref{p7} also for $t=\infty$.
This concludes the proof.
\end{proof}

For any $t\in[0,\infty]$ let $\tQQ_t$ denote the absolutely continuous part
of the measure $\tPP_t$ with respect to~$\PP_t$. Let us note that
$\tQQ_t(\Omega)=1\Longleftrightarrow\tPP_t\ll\PP_t$.

\begin{lemma}
\label{PP4}
We have
\begin{equation*}
\frac{d\tQQ_t}{d\PP_t}=Z_t\quad\PP\text{-a.s.},\quad t\in[0,\infty].
\end{equation*}
\end{lemma}

\begin{proof}
Let us fix $t\in[0,\infty]$.
By Lemma~\ref{PP1} and Jessen's theorem (see \cite[Th.~5.2.26]{Stroock:93}),
\begin{equation*}
\frac{d\tPP_{t\wedge\tau_n}}{d\PP_{t\wedge\tau_n}}\xrightarrow[n\uparrow\infty]{}
\frac{d\tQQ_t}{d\PP_t}\quad\PP\text{-a.s.}
\end{equation*}
Since the process $Z$ is stopped at~$\zeta$,
we have $Z_{t\wedge\tau_n}\to Z_t$ $\PP$-a.s. as~$n\uparrow\infty$.
Now the statement follows from Lemma~\ref{PP3}.
\end{proof}

\begin{proof}[Proof of Theorem~\ref{MRA1}]
We prove the remark after Theorem~\ref{MRA1}
(from this remark Theorem~\ref{MRA1} follows).
By Lemma~\ref{PP4} we have
\begin{equation*}
(Z_t)_{t\in[0,T]}\text{ is a }\PP\text{-martingale}
\Longleftrightarrow
\EE_\PP\,\frac{d\tQQ_T}{d\PP_T}=1
\Longleftrightarrow
\tPP_T\ll\PP_T
\Longleftrightarrow
S>T\;\;\tPP\text{-a.s.},
\end{equation*}
where the last equivalence follows from the remark after Lemma~\ref{ST4}.
Assume first that $\PP\ne\tPP$ which is, by the occupation times formula,
equivalent to $\nu_L(b\ne0)>0$.
Let us recall that all points in $J$ are non-separating in our case
and that an endpoint of $J$ is non-separating if and only if it is good
in the sense of \eqref{mra12} and~\eqref{mra14}.
Applying now remark~(ii) following Theorem~\ref{A0} we get
that $S>T$~$\tPP$-a.s. if and only if the coordinate process $X$
does not exit~$J$ at a bad endpoint under $\tPP$. 
This gives the criterion in Theorem~\ref{MRA1}.

Finally, let us assume that $\nu_L(b\ne0)=0$ (i.e. $\PP=\tPP$).
Then, clearly, if $l$ (resp.~$r$) is bad, then $\wt s(l)=-\infty$ (resp. $\wt s(r)=\infty$),
hence $X$ does not exit at~$l$ (resp. at~$r$) under~$\tPP$.
This means that the criterion in Theorem~\ref{MRA1}
works also in the case $\PP=\tPP$. This proves the theorem.
\end{proof}

\begin{proof}[Proof of Theorem~\ref{MRA2}]
By Lemma~\ref{PP4} we have
\begin{equation*}
(Z_t)_{t\in[0,\infty)}\text{ is a u.i. }\PP\text{-martingale}
\Longleftrightarrow
\EE_\PP\,\frac{d\tQQ_\infty}{d\PP_\infty}=1
\Longleftrightarrow
\tPP\ll\PP
\Longleftrightarrow
S=\delta\;\;\tPP\text{-a.s.},
\end{equation*}
where 
``u.i.'' stands for ``uniformly integrable''
and
the last equivalence follows from Lemma~\ref{ST4}.
We can now apply Theorem~\ref{A0} and remark~(ii) after it,
in much the same way as in the proof of Theorem~\ref{MRA1},
to express the condition $S=\delta$~$\tPP$-a.s.
in terms of the functions $\mu$, $\sigma$, and~$b$.
In the case $\PP\ne\tPP$ (i.e. $\nu_L(b\ne0)>0$),
we find that $S=\delta$~$\tPP$-a.s.
if and only if either both endpoints of $J$ are good
or $\tPP$-almost all paths of $X$ converge to a good endpoint of $J$ as $t\uparrow\zeta$.
In the latter possibility the other endpoint of $J$ may be (and, actually, will be) bad.
The former possibility is condition~(D) in Theorem~\ref{MRA2}.
It follows from Proposition~\ref{APP2} that the latter possibility
is described by conditions (B) and~(C).
Finally, condition~(A) constitutes the remaining case $\nu_L(b\ne0)=0$.
This concludes the proof.
\end{proof}

\appendix
\section{Local Time of One-Dimensional Diffusions}
\label{sec:a:lt}
In this appendix we describe some properties of the local time of solutions
of one-dimensional SDEs with the coefficients satisfying the Engelbert--Schmidt conditions.
These properties are used in Section~\ref{MR}.
Consider 
a $J$-valued diffusion $Y=(Y_t)_{t\in[0,\infty)}$,
where
$J=(l,r)$, $-\infty\le l<r\le\infty$,
governed by SDE~\eqref{mr1}
on some filtered probability space
$(\Omega,\cF,(\cF_t)_{t\in[0,\infty)},\PP)$,
where the coefficients $\mu$ and $\sigma$
are Borel functions $J\to\bbR$ satisfying \eqref{mr2}--\eqref{mr3}.
Since $Y$ is a continuous semimartingale up to the exit time
$\zeta$, one can define local time $\{L_t^y(Y)\colon y\in J,t\in[0,\zeta)\}$
on the stochastic interval $[0,\zeta)$ for any $y\in J$
in the usual way (e.g. via the obvious generalisation of~\cite[Ch.~VI, Th.~1.2]{RevuzYor:99}).

It follows from Theorem~VI.1.7 in~\cite{RevuzYor:99}
that the random field $\{L_t^y(Y)\colon y\in J,t\in[0,\zeta)\}$
admits a modification such that the map $(y,t)\mapsto L^y_t(Y)$
is a.s. continuous in~$t$ and c\`adl\`ag in~$y$.
As usual we always work with this modification.
First, and this is also of interest in itself,
we show that for solutions of SDEs with the coefficients
satisfying the Engelbert--Schmidt conditions,
the local time is moreover jointly continuous in $t$ and~$y$.
(Let us note that joint continuity is a stronger property
than continuity in $t$ for any $y$ and continuity in $y$ for any~$t$.)

\begin{proposition}
\label{LocTime_cont}
The random field
$\{L_t^y(Y)\colon y\in J,t\in[0,\zeta)\}$
has a modification such that the map
$$
J\times[0,\zeta)\ni(y,t)\mapsto L_t^y(Y)\in\bbR_+
$$
is a.s. jointly continuous in $(t,y)$.
Furthermore, any modification of the local time that is a.s. continuous in~$t$ and c\`adl\`ag in~$y$
is, in fact, a.s. jointly continuous in $(t,y)$.
\end{proposition}

\begin{proof}
Tanaka's formula~\cite[Ch.~VI, Th.~1.2]{RevuzYor:99}
implies for any $y\in J$ that
\begin{equation}
\label{eq:a:lt1}
L_t^y(Y)=2\left[(Y_t-y)^+-(Y_0-y)^+-\int_0^tI_{\{Y_s>y\}}\sigma(Y_s)\dd W_s-\int_0^t I_{\{Y_s>y\}}\mu(Y_s)\dd s\right]
\end{equation}
a.s. on $\{t<\zeta\}$, where $I$ denotes the indicator function.
Thus, we need to prove that the right-hand side of~\eqref{eq:a:lt1}
admits an a.s. jointly continuous in $t$ and $y$ modification.
This is clear for the first two terms there.
It is shown in the proof of Theorem VI.1.7 of~\cite{RevuzYor:99}
that Kolmogorov's criterion for the existence of continuous modifications of
random fields (see~\cite[Ch.~II, Th.~2.1]{RevuzYor:99})
applies to the integral with respect to $W$ in the right-hand side of~\eqref{eq:a:lt1}.
It remains to deal with the term $\int_0^t I_{\{Y_s>y\}}\mu(Y_s)\,ds$,
which, by~\eqref{mr2} and the occupation times formula, is equal to
$$
F(y,t):=\int_y^r \frac\mu{\sigma^2}(x)L^x_t(Y)\,dx
$$
a.s. on $\{t<\zeta\}$.
Note that the latter expression is finite due to~\eqref{mr3}
and the fact that the map $x\mapsto L^x_t(Y)$ has a.s. on $\{t<\zeta\}$
a compact support in $J$ (see \cite[Ch.~VI, Prop.~(1.3)]{RevuzYor:99})
and is c\`adl\`ag, hence bounded.
It suffices now to show that for any $\eps>0$,
$F$ is jointly continuous a.s. on $\{t+\eps<\zeta\}$.
Let us take a small $\delta>0$.
Then for $(y',t')\in J\times[0,\zeta)$
such that $|y-y'|<\delta$ and $|t-t'|<\eps$,
we have
\begin{eqnarray*}
\lvert F(y,t)-F(y',t')\rvert
& \leq &
\lvert F(y,t)-F(y',t)\rvert+
\lvert F(y',t)-F(y',t')\rvert\\
& \leq & \int_{y-\delta}^{y+\delta}\frac{|\mu|}{\sigma^2}(x) L^x_t(Y)\dd x + 
\int_l^r\frac{|\mu|}{\sigma^2}(x) \lvert L_t^x(Y)-L_{t'}^x(Y)\rvert\dd x.
\end{eqnarray*}
Again, due to~\eqref{mr3} and boundedness of the mapping $x\mapsto L^x_t(Y)$,
the first integral is arbitrarily small provided $\delta$ is small enough.
The second integral tends to~$0$ as $t'\to t$ by the dominated convergence theorem
(as above we have $\int_l^r (|\mu|/\sigma^2)(x)L^x_{t+\eps}(Y)\,dx<\infty$ a.s. on $\{t+\eps<\zeta\}$).
We proved that the local time has an a.s. jointly continuous modification.

To prove the last statement let us take a version $\{L_t^y(Y)\colon y\in J,t\in[0,\zeta)\}$
of the local time, which is a.s. continuous in $t$ and c\`adl\`ag in~$y$,
and consider its modification $\{\wt L_t^y(Y)\colon y\in J,t\in[0,\zeta)\}$,
which a.s. is jointly continuous in $(t,y)$. Then a.s. we have
\begin{equation}
\label{eq:a:lt2}
L^y_t(Y)=\wt L^y_t(Y)\text{ for all }y\in J\cap\bbQ,\;t\in\bbQ_+,\;t<\zeta.
\end{equation}
Since both $(L^y_t(Y))$ and $(\wt L^y_t(Y))$ are continuous in $t$ and c\`adl\`ag in~$y$,
\eqref{eq:a:lt2}~implies that they are indistinguishable.
Hence, $(L^y_t(Y))$ is, in fact, a.s. jointly continuous in $t$ and~$y$.
\end{proof}

The next result is a consequence of Theorem~2.7 in~\cite{ChernyEngelbert:05}.

\begin{proposition}
\label{prop:a:lt2}
Let $\alpha\in J$ and set
$$
\tau_\alpha:=\inf\{t\in[0,\infty)\colon Y_t=\alpha\}
$$
(with $\inf\emptyset:=\infty$). Then $L^\alpha_t(Y)>0$ a.s. on $\{\tau_\alpha<t<\zeta\}$.
\end{proposition}

\begin{remark}
By Proposition~\ref{APP0}, $\PP(\tau_\alpha<\zeta)>0$.
Hence, there exists $t\in(0,\infty)$ such that $\PP(\tau_\alpha<t<\zeta)>0$.
\end{remark}

Let us note that the result of Proposition~\ref{prop:a:lt2} no longer holds
if the coefficients $\mu$ and $\sigma$ of~\eqref{mr1} fail to satisfy the Engelbert--Schmidt conditions
(see Theorem~2.6 in~\cite{ChernyEngelbert:05}).

\section{Behaviour of One-Dimensional Diffusions}
\label{sec:a:sde}
Here we state some well-known results
about the behaviour of solutions of one-dimensional SDEs
with the coefficients satisfying the Engelbert--Schmidt conditions
(see e.g.~\cite[Sec.~4.1]{ChernyEngelbert:05})
that are extensively used in Sections \ref{STB} and~\ref{P}.

Let $(\Omega,\cF,(\cF_t)_{t\in[0,\infty)})$
be the canonical filtered space described
in Section~\ref{STB}, $X$ the coordinate process on~$\Omega$,
and $\zeta$ the exit time of~$X$.
Let the measure $\PP$ on $(\Omega,\cF)$ be the distribution of a (unique in law)
solution of SDE~\eqref{mr1}, where the coefficients
$\mu,\sigma$ are Borel functions $J\to\bbR$ satisfying
conditions \eqref{mr2}--\eqref{mr3}.

Consider the sets
\begin{align*}
A&=\bigl\{\,\zeta=\infty,\;
\limsup_{t\to\infty} X_t=r,\;
\liminf_{t\to\infty} X_t=l\bigr\},\\
B_r&=\bigl\{\,\zeta=\infty,\;
\lim_{t\to\infty} X_t=r\bigr\},\\
C_r&=\bigl\{\,\zeta<\infty,\;
\lim_{t\uparrow\zeta} X_t=r\bigr\},\\
B_l&=\bigl\{\,\zeta=\infty,\;
\lim_{t\to\infty} X_t=l\bigr\},\\
C_l&=\bigl\{\,\zeta<\infty,\;
\lim_{t\uparrow\zeta} X_t=l\bigr\}.
\end{align*}
For $a\in J$ define the stopping time
\begin{equation*}
\tau_a=\inf\{t\in[0,\infty)\colon X_t=a\}\quad(\inf\emptyset:=\infty).
\end{equation*}

\begin{proposition}
\label{APP0}
For any $a\in J$ we have $\PP(\tau_a<\infty)>0$.
\end{proposition}

\begin{proposition}
\label{APP1}
Either $\PP(A)=1$ or
$\PP(B_r\cup B_l\cup C_r\cup C_l)=1$.
\end{proposition}

Let $s$ denote the scale function of diffusion~\eqref{mr1}
and $\rho$ the derivative of~$s$ (see \eqref{mra8} and~\eqref{mra9}).

\begin{proposition}
\label{APP2}
(i) If $s(r)=\infty$, then
$\PP(B_r\cup C_r)=0$.

(ii) Assume that $s(r)<\infty$.
Then either $\PP(B_r)>0$, $\PP(C_r)=0$
or $\PP(B_r)=0$, $\PP(C_r)>0$.
Moreover, we have
\begin{equation*}
\PP\bigl(\,\lim_{t\uparrow\zeta}X_t=r,\;X_t>a\;\forall t\in[0,\zeta)\bigr)>0
\end{equation*}
for any $a<x_0$.
\end{proposition}

Condition~\eqref{mra21}
is necessary and sufficient for $\PP(B_r)=0$, $\PP(C_r)>0$.
Proposition~\ref{APP2}
therefore implies that a necessary and sufficient condition for 
$\PP(B_r)>0$, $\PP(C_r)=0$
consists of $s(r)<\infty$ and
$\frac{s(r)-s}{\rho\sigma^2}\notin\Llocr$.

Proposition~\ref{APP2} concerns the behaviour
of the solution of SDE~\eqref{mr1} at the endpoint~$r$.
Clearly, it has its analogue for the behaviour at~$l$.

\bibliographystyle{abbrv}
\bibliography{refs}
\end{document}